\documentclass[11pt, hidelinks]{article} 

\usepackage{hyperref, mathscinet}
\usepackage{url}
\usepackage{fullpage}
\usepackage{xcolor}
\usepackage{enumitem, comment, xifthen}
\usepackage{graphicx}
\graphicspath{{figs/}}
\usepackage{subfig}
\usepackage{overpic}
\usepackage{booktabs,float}
\usepackage{amsmath,amssymb,amsthm, amsfonts}
\usepackage[T1]{fontenc}
\theoremstyle{plain}

\newtheorem{theo}{Theorem}

\newtheorem{lem}{Lemma}
\newtheorem{prop}{Proposition}
\newtheorem{cor}{Corollary}

\theoremstyle{definition} 

\newtheorem{nota}{Notation}
\newtheorem{de}{Definition}
\newtheorem{exa}{Example}
\newtheorem{as}{Assumption}
\newtheorem{alg}{Algorithm}

\newcommand{\btheo}{\begin{theo}}
\newcommand{\bde}{\begin{de}}
\newcommand{\ble}{\begin{lem}}
\newcommand{\bpr}{\begin{prop}}
\newcommand{\bno}{\begin{nota}}
\newcommand{\bex}{\begin{exa}}
\newcommand{\bcor}{\begin{cor}}
\newcommand{\spro}{\begin{proof}}
\newcommand{\bas}{\begin{as}}
\newcommand{\balg}{\begin{alg}}

\newcommand{\etheo}{\end{theo}}
\newcommand{\ede}{\end{de}}
\newcommand{\ele}{\end{lem}}
\newcommand{\epr}{\end{prop}}
\newcommand{\eno}{\end{nota}}
\newcommand{\eex}{\end{exa}}
\newcommand{\ecor}{\end{cor}}
\newcommand{\fpro}{\end{proof}}
\newcommand{\eas}{\end{as}}
\newcommand{\ealg}{\end{alg}}

\theoremstyle{plain}

\newtheorem{theos}{Theorem}
\newtheorem{props}{Proposition}
\newtheorem{lems}{Lemma}
\newtheorem{cors}{Corollary}

\theoremstyle{definition}
\newtheorem{exas}{Example}
\newtheorem{algs}{Algorithm}
\newtheorem{asss}{Assumption}
\newtheorem{defns}{Definition}

\newcommand{\btheos}{\begin{theos}}
\newcommand{\etheos}{\end{theos}}
\newcommand{\bprops}{\begin{props}}
\newcommand{\eprops}{\end{props}}
\newcommand{\bdes}{\begin{defns}}
\newcommand{\edes}{\end{defns}}
\newcommand{\blems}{\begin{lems}}
\newcommand{\elems}{\end{lems}}
\newcommand{\bcors}{\begin{cors}}
\newcommand{\ecors}{\end{cors}}
\newcommand{\bexs}{\begin{exas}}
\newcommand{\eexs}{\end{exas}}
\newcommand{\balgs}{\begin{algs}}
\newcommand{\ealgs}{\end{algs}}
\newcommand{\bass}{\begin{asss}}
\newcommand{\eass}{\end{asss}}

\usepackage{hyperref}
\usepackage{url}
\usepackage[margin=1in]{geometry}
\usepackage{amsmath,amssymb,amsthm, amsfonts}
\usepackage[T1]{fontenc}
\usepackage{mathtools}
\usepackage{xcolor}
\usepackage{enumitem, comment, xifthen}
\usepackage{graphicx}
\graphicspath{{figs/}}
\usepackage{subfig}
\theoremstyle{plain}  
 \newtheorem{definition}{\textbf{Definition}}
 \theoremstyle{definition}

\theoremstyle{remark}

\makeatletter
\def\letterdef#1#2#3{\def\letterdef@##1{\expandafter\def\csname #1\endcsname{#2}}%
  \letterdef@@#3{?\@car{}}\@nil}
\def\letterdef@@#1{\@gobble#1\letterdef@{#1}\letterdef@@}
\makeatother

\makeatletter
\newcommand*{\coloneq}{\mathrel{\rlap{%
                     \raisebox{0.3ex}{$\m@th\cdot$}}%
                     \raisebox{-0.3ex}{$\m@th\cdot$}}%
  =}
\makeatother
\makeatletter
\newcommand*{\eqcolon}{=\mathrel{\rlap{%
                     \raisebox{0.3ex}{$\m@th\cdot$}}%
                     \raisebox{-0.3ex}{$\m@th\cdot$}}%
}
\makeatother
\DeclarePairedDelimiterX{\klx}[2]{(}{)}{%
  #1\;\delimsize\|\;#2%
}
\DeclarePairedDelimiterX{\quantklx}[3]{(}{)}{%
  #1\;\delimsize\|\;#2\;\delimsize\vert\;#3%
}
\DeclarePairedDelimiterX{\inner}[2]{\langle}{\rangle}{%
  #1,#2%
}

\newcommand{\R}{\mathbf R} 

\letterdef{c#1}{\mathcal{#1}}{ABCDEFGHIJKLMNOPQRSTUVWXYZ}
\let\defn\coloneq

\newcommand{\ceil}[1]{\left\lceil #1 \right\rceil}

\newcommand{\ud}[0]{\,\mathrm{d}}  
\newcommand{\1}{\mathbf 1} 
\let\ones\1
 
\let\epsilon\varepsilon
\newcommand{\eps}{\varepsilon}

\let\hat\widehat

\renewcommand{\leq}{\leqslant}
\renewcommand{\geq}{\geqslant}

\newcommand{\kl}{D_{\rm kl}\klx}

\newcommand{\iid}{\textnormal{i.i.d.}}

\newcommand{\simiid}{\stackrel{\iid}{\sim}} 

\makeatletter
\newcommand{\E}{\operatorname*{\mathbf{E}}\ilimits@}
\makeatother
\makeatletter
\renewcommand{\P}{\operatorname*{\mathbf{P}}\ilimits@}
\makeatother

\newcommand\blfootnote[1]{%
  \begingroup
  \renewcommand\thefootnote{}\footnote{#1}%
  \addtocounter{footnote}{-1}%
  \endgroup
}

\newcounter{algorithmctr}
\renewcommand{\thealgorithmctr}{\arabic{algorithmctr}}
   {\refstepcounter{algorithmctr}\begin{list}{}{%
       \setlength{\rightmargin}{0\linewidth}%
       \setlength{\leftmargin}{0\linewidth}}%
       \rmfamily\small
       \item[]{\setlength{\parskip}{0ex}\hrulefill\par%
        \nopagebreak{\bfseries\textsf{Algorithm \thealgorithmctr~}}}}%
   {{\setlength{\parskip}{-1ex}\nopagebreak\par\hrulefill} \end{list}}

\makeatletter
\long\def\@makecaption#1#2{
        \vskip 0.8ex
        \setbox\@tempboxa\hbox{\small {\bf #1.} #2}
        \parindent 1.5em 
        \dimen0=\hsize
        \advance\dimen0 by -3em
        \ifdim \wd\@tempboxa >\dimen0
                \hbox to \hsize{
                        \parindent 0em
                        \hfil 
                        \parbox{\dimen0}{\def\baselinestretch{0.96}\small
                                {\bf #1.} #2
                                } 
                        \hfil}
        \else \hbox to \hsize{\hfil \box\@tempboxa \hfil}
        \fi
        }
\makeatother

\newcommand{\upstairs}[1]{\textsuperscript{#1}}
\newcommand{\affilone}{\dag}
\newcommand{\affiltwo}{\ddag}

\newcommand{\affilchic}{$\diamond$}
\newcommand{\fstar}{f^\star}
\newcommand{\xspace}{\ensuremath{\cX}}
\newcommand{\noise}{\ensuremath{\xi}}
\newcommand{\noisebound}{\sigma}
\newcommand{\sourcedist}{P}
\newcommand{\targetdist}{Q}
\newcommand{\dimension}{k} 
\newcommand{\numobs}{n}

\newcommand{\numsource}{{\numobs_\sourcedist}}
\newcommand{\numtarget}{{\numobs_\targetdist}}
\newcommand{\ball}[2]{\mathsf{B}(#1, #2)}
\newcommand{\bandwidth}{h}
\newcommand{\similarity}[3]{\rho_{#1}(#2 , #3)}
\newcommand{\simmeasgeneric}{\rho_\bandwidth}
\newcommand{\fhat}{\ensuremath{\widehat{f}}}
\newcommand{\fhatNW}{\ensuremath{\fhat}}
\newcommand{\fbar}{\ensuremath{\overline{f}}}
\newcommand{\goodset}[1]{\ensuremath{\cG_{#1}}}

\newcommand{\midpointdist}{\mu}
\newcommand{\upperconst}{c_u}
\newcommand{\lowerconst}{c_\ell}
\newcommand{\caseindicator}{\eta}
\newcommand{\holderclass}{\cF(\beta, L)}
\newcommand{\lowerboundclassbig}[2]{\cD(#1, #2)}
\newcommand{\lowerboundclasssmall}[2]{\cD'(#1, #2)}
\newcommand{\supportlength}{S}
\newcommand{\Bin}[2]{\mathsf{Bin}(#1, #2)}
\newcommand{\Normal}[2]{\mathsf{N}(#1, #2)}

\newcommand{\randbin}{b}
\newcommand{\GVset}{\cB}
\newcommand{\psiintegral}{C_\Psi}
\newcommand{\numupper}{n_u}
\newcommand{\numlower}{n_\ell}
\newcommand{\numeff}{\numobs_\mathrm{eff}}

\newcommand{\transClass}{\mathcal{T}}

\renewcommand\blfootnote[1]{%
  \begingroup
  \renewcommand\thefootnote{}\footnote{#1}%
  \addtocounter{footnote}{-1}%
  \endgroup
}

\newcommand{\trad}{\ensuremath{K}}
\newcommand{\polyexp}{\ensuremath{\kappa}}
\newcommand{\Term}{\ensuremath{T}}
\newcommand{\jdist}{\ensuremath{\nu}}

\begin{document}
\begin{center}
  {\bf{\LARGE A new similarity measure for covariate shift \\
      with applications to nonparametric regression}}\\
  
  \vspace*{.2in}

  \begin{tabular}{cc}
    Reese Pathak\upstairs{\affilone,\! $\star$}, 
    Cong Ma\upstairs{\affilchic,\! $\star$}, 
    and
    Martin J.\ Wainwright\upstairs{\affilone, \affiltwo} \\[1.5ex]
    \upstairs{\affilchic} Department of Statistics, University of Chicago \\
    \upstairs{\affilone} Department of Electrical Engineering and Computer
    Sciences, UC Berkeley \\
    \upstairs{\affiltwo} Department of Statistics, UC Berkeley \\
    \texttt{congm@uchicago.edu, \string{pathakr,wainwrig\string}@berkeley.edu}
  \end{tabular}
  \vspace*{.2in}
  \begin{abstract}
    We study covariate shift in the context of nonparametric
    regression.  We introduce a new measure of distribution mismatch
    between the source and target distributions that is based on the
    integrated ratio of probabilities of balls at a given radius.  We
    use the scaling of this measure with respect to the radius to
    characterize the minimax rate of estimation over a family of
    Hölder continuous functions under covariate shift.  In comparison
    to the recently proposed notion of transfer exponent, this measure
    leads to a sharper rate of convergence and is more
    fine-grained. We accompany our theory with concrete instances of
    covariate shift that illustrate this sharp difference.
  \end{abstract}
\end{center}


\section{Introduction}

In\blfootnote{\upstairs{$\star$} RP and CM contributed equally to this
  paper.}  the standard formulation of prediction or classification,
future data (as represented by a test set) is assumed to be drawn
from the same distribution as the training data.  This assumption,
while theoretically convenient, may fail to hold in many real-world
scenarios.  For instance, training data might be collected only from a
sub-group within a broader population (such as in medical trials), or
the environment might change over time as data are collected.  Such
scenarios result in a distribution mismatch between the training and
test data.

In this paper, we study an important case of such distribution
mismatch---namely, the covariate shift problem (e.g.,~\cite{Shi00,
  QuiEtAl09}).  Suppose that a statistician observes
covariate-response pairs $(X, Y)$, and wishes to build a prediction
rule.  In the problem of covariate shift, the distribution of the
covariates $X$ is allowed to change between the training and test
data, while the posterior distribution of the responses (namely, $Y
\mid X$) remains fixed.  Compared to the usual i.i.d.\ setting, this
serves as a more accurate model for a variety of real-world
applications, including image classification~\cite{SaeEtAl2010},
biomedical engineering~\cite{LiEtAl2010}, sentiment
analysis~\cite{BliEtAl07}, and audio processing~\cite{HasEtAl13},
among many others.

More formally, suppose that the statistician observes $\numsource$
covariates $\{X_i\}_{i=1}^{\numsource}$ from a \emph{source
distribution} $\sourcedist$, and $\numtarget$ covariates
$\{X_i\}_{i=\numsource + 1}^{\numtarget + \numsource}$ from a
\emph{target distribution} $\targetdist$. For each observed $X_i$, she
also observes a response $Y_i$ drawn from the same conditional
distribution.  The \emph{regression function} $\fstar(x) = \E[Y \mid
  x]$ defined by this conditional distribution is assumed to lie in
some function class $\cF$.  The statistician uses these samples to
produce an estimate $\fhat$, which will be evaluated on the target
distribution, with a fresh sample $X \sim \targetdist$, yielding
the mean-squared error
\begin{equation*}
\|\hat f - \fstar\|_{L^2(Q)}^2 \defn \E\Big[\big(\hat f(X) -
  \fstar(X)\big)^2\Big].
\end{equation*}
When there is no covariate shift, the fundamental (minimax) risks for
this problem are well-understood~\cite{Hal72, IbrKha80, Sto82}.  The
goal of this paper is to understand how, for nonparametric function
classes $\cF$, this minimax risk changes as a function of the
``amount'' of covariate shift between $\sourcedist$ and $\targetdist$.

\subsection{Our contributions and related work}

Let us summarize the main contributions of this paper, and put them in
the context of related work.

\paragraph{Our contributions.} We introduce a similarity measure\footnote{To be clear, this quantity  actually serves as a \emph{dis}-similarity measure: as shown in the sequel, source-target pairs $(P, Q)$
with larger values $\similarity{\bandwidth}{\sourcedist}{\targetdist}$
lead to ``harder'' estimation problems in terms of covariate shift.}
$\simmeasgeneric$ between two probability measures $\sourcedist,
\targetdist$ on a common metric space $(\xspace, d)$. For any level
$\bandwidth > 0$, it is defined as
\begin{equation}
\label{eqn:similarity-measure}
\similarity{\bandwidth}{\sourcedist}{\targetdist} \defn \int_\xspace
\frac{1}{\sourcedist\big(\ball{x}{\bandwidth}\big)} \, \ud
\targetdist(x),
\end{equation}
where $\ball{x}{\bandwidth} \defn \{\, x' \in \xspace \mid d(x, x')
\leq \bandwidth\,\}$ is the closed ball of radius $\bandwidth$
centered around $x$.  We substantiate the significance of this
similarity measure via the following contributions:
\begin{enumerate}[label=
{\upshape(\roman*)}]
\item For regression functions that are H\"{o}lder continuous, we
  demonstrate a performance guarantee for the Nadaraya-Watson
  kernel estimator under covariate shift that is fully determined by the
  scaling of the similarity measure
  $\similarity{\bandwidth}{\sourcedist}{\targetdist}$ with respect to
  the radius $\bandwidth$.
\item We complement these upper bounds with matching lower bounds---in
  a minimax sense---demonstrating that the best achievable rate of
  estimation in Hölder classes is also determined by the scaling of
  this similarity measure.
\item We show how the similarity measure $\simmeasgeneric$ can be
  controlled based on the metric properties of the space $\xspace$. In
  addition, we compare $\simmeasgeneric$ with existing notions for
  covariate shift (e.g., bounded likelihood ratios, transfer
  exponents), thereby showcasing some of its advantages.
\end{enumerate}

\paragraph{Related work.} 
The problem of covariate shift was studied in the seminal work by
Shimodaira~\cite{Shi00}, who provided asymptotic guarantees for a
weighted maximum likelihood estimator under covariate shift.  Since
then, a plethora of work has analyzed covariate shift, or the general
distribution mismatch problem (also referred to as domain adaptation
or transfer learning).

For general distribution mismatch, one line of work provides rates
that depend on distance metrics between the source-target pair (e.g.,
~\cite{BenBliEtAl10,BenLuEtAl10,
  GerEtAl13,ManEtAl09,CorEtAl19,MohEtAl12}). These results hold under
fairly general conditions, but do not necessarily guarantee
consistency as the sample size $n$ increases. In contrast, our
guarantees for covariate shift do guarantee consistency, and moreover,
we provide explicit nonasymptotic, optimal nonparametric rates.  As
pointed out in the paper~\cite{KpoMar21}, the distribution mismatch
problem is asymmetric in the sense that it may be easier to estimate
accurately when dealing with covariate shift from $\sourcedist$ to
$\targetdist$ than from $\targetdist$ to $\sourcedist$. Our results
also corroborate this intuition.  It is worth noting that these prior
distance metrics fall short of capturing the inherent asymmetry
between $\sourcedist$ and $\targetdist$.

Another line of work addresses covariate shift under conditions on the
likelihood ratio $\ud\targetdist /\ud\sourcedist$.  For instance, some
authors have obtained results for bounded likelihood
ratios~\cite{SugEtAl12,Kpo17} or in terms of information-theoretic
divergences between the source-target
pair~\cite{SugNakEtAl08,ManMohEtAl09}.  Our work is inspired in part
by the work of Kpotufe and Martinet~\cite{KpoMar21}, who introduced
the notion of the \emph{transfer exponent}.  It is a condition that
bounds the mass placed by the pair $(\sourcedist, \targetdist)$ on
balls of varying radii; using this notion, they analyzed various
problems of nonparametric classification.  Our work, focusing instead
on nonparametric regression problems and using the measure
$\simmeasgeneric$, provides sharper rates than those obtainable by
considering the transfer exponent; see Section~\ref{sec:comparison}
for details. Thus, the similarity measure $\simmeasgeneric$ provides a
more fine-grained control on the effect of covariate shift on
nonparametric regression.

Finally, it is worth mentioning other recent works that give risk
bounds for covariate shift problems, including on linear
models~\cite{LeiHuLee21}, as well as linear models and one-layer
neural networks~\cite{MouMohEtAl20}.  Although these results deal with
covariate shift, the rates obtained are parametric ones, and hence not
directly comparable to the nonparametric rates that are the focus of
our inquiry.


\subsection{Notation}
Here we collect notation used throughout the paper.  We use $\R$ to
denote the real numbers.  We use $(\xspace, d)$ to denote a metric
space, and we equip it with the usual Borel $\sigma$-algebra.  We let
\mbox{$\ball{x}{r} \defn \big \{\, x' \in \xspace \mid d(x, x') \leq r
  \, \big \}$} be the closed ball of radius $r$ centered at $x$.  We
reserve the capital letters $X, Y$, possibly with subscripts, for a
pair of random variables arising from a regression model.  Similarly,
we reserve $\sourcedist, \targetdist$ for a pair of two probability
measures on $(\xspace, d)$.  For $\bandwidth > 0$, we denote by
$N(\bandwidth)$ the covering number of $\xspace$ at resolution
$\bandwidth$ in the metric $d$.  This is the minimal number of balls
of radius at most $\bandwidth > 0$ required to cover the space
$\xspace$. \\

The remainder of this paper is organized as follows.  We begin in
Section~\ref{SecMain} by setting up the problem more precisely, and
stating and discussing our main results on covariate shift: namely,
upper bounds in Theorem~\ref{thm:upper-bound}, accompanied by matching
lower bounds in Theorem~\ref{thm:lower-bound}.  These results
establish that the similarity measure~\eqref{eqn:similarity-measure}
provides a useful measure of the ``difficulty'' of source-target pairs
in covariate shift; accordingly, Section~\ref{SecProperties} is
devoted to a comparison and discussion of this measure relevant to
concepts from past work, including likelihood ratio bounds and
transfer exponents.  The proofs of all our results are given in
Section~\ref{SecProofs}, and we conclude with a discussion in
Section~\ref{SecDiscussion}.


\section{How covariate shift affects nonparametric regression}
\label{SecMain}

In this section, we use the similarity measure introduced in
equation~\eqref{eqn:similarity-measure} to characterize how covariate
shift can change the minimax risks of estimation for certain classes
of nonparametric regression models.  We begin in
Section~\ref{SecBackground} by setting up the observation model to be
considered, along with some associated assumptions on the regression
function $\fstar$, the conditional distribution of $Y \mid X$, and the
covariate shift.  In Section~\ref{SecNW}, we derive an achievable
result (Theorem~\ref{thm:upper-bound}) for nonparametric regression in
the presence of covariate shift, in particular via a careful analysis
of the classical Nadaraya-Watson estimator.  Our upper bound in this
section is general, and illustrates the key role of the similarity
measure $\simmeasgeneric$.  In Section~\ref{SecAlphaFamily}, we
introduce the $\alpha$-families of source-target pairs $(P, Q)$, and
use Theorem~\ref{thm:upper-bound} to derive achievable results for
these families.  In Section~\ref{SecLower}, we state some
complementary lower bounds for $\alpha$-families
(Theorem~\ref{thm:lower-bound}), showing that our achievable results
are, in fact, unimprovable.


\subsection{Observation  model and assumptions}
\label{SecBackground}
Suppose that we observe covariate-response pairs $\{(X_i,
Y_i)\}_{i=1}^\numobs \subset \xspace \times \R$ that are drawn from
nonparametric regression model of the following type.  The conditional
distribution of $Y \mid X$ is the same for all $i = 1, \dots,
\numobs$, and our goal is to estimate the regression function
$\fstar(x) \defn \E[Y \mid X = x]$.  In terms of the ``noise''
variables, $\noise_i \defn Y_i - f^\star(X_i)$, the observations can
be written in the form
\begin{equation}
\label{EqnObsModel}  
Y_i = \fstar(X_i) + \noise_i, \quad i = 1, \ldots, \numobs.
\end{equation}
In our analysis, we impose three types of regularity conditions: (i)
H\"{o}lder continuity of the regression function; (ii) the type of
covariate shift allowed; and (iii) tail conditions on the noise
variables $\{\noise_i\}_{i=1}^\numobs$.
\bas[H\"{o}lder continuity]
\label{asmp:holder}
For some $L > 0$ and $\beta \in (0, 1]$, the function $\fstar
  \colon \xspace \to \R$ is $(\beta, L)$-H\"{o}lder continuous,
  meaning that
\begin{equation*}
\big|\fstar(z) - \fstar(z')\big| \leq L \; [d(z, z')]^\beta, \quad
\mbox{for any}~z, z' \in \xspace.
\end{equation*}
\eas
\noindent We note that in the special case $\beta = 1$, the function
$\fstar$ is $L$-Lipschitz.

\bas[Covariate shift]
\label{asmp:covshift}
The covariates $X_1, \dots, X_\numobs$ are independent, and drawn as
\begin{equation*}
X_1, \dots, X_{\numsource} \simiid \sourcedist \quad \mbox{and} \quad
X_{\numsource + 1}, \dots X_{\numsource + \numtarget} \simiid
\targetdist \qquad \mbox{where $\numobs = \numsource + \numtarget$.}
\end{equation*}
\eas

\bas[Noise assumption]
\label{asmp:noise}
The variables $\{\noise_i\}_{i=1}^\numobs$ satisfy the second moment
bound
\begin{equation*}
  \sup_x \E\big[\noise_i^2 \mid X_i = x\big] \leq \noisebound^2 \qquad
  \mbox{for $i = 1, \ldots, \numobs$.}
\end{equation*}
\eas
\noindent Note that by construction, the variables $\noise_i$ are
(conditionally) centered.  Assumption~\ref{asmp:noise} also allows
$\noise_i$ to depend on $X_i$, as long as the variance is uniformly
bounded above.

\subsection{Achievable performance via the Nadaraya-Watson estimator}
\label{SecNW}

We first exhibit an achievable result for the problem of nonparametric
regression in the presence of covariate shift.  We do so by analyzing
a classical and simple method for nonparametric estimation, namely the
Nadaraya-Watson estimator~\cite{Nad64, Wat64}, or NW for short. The
main result of this section is to show that the mean-squared error
(MSE) of the NW estimator is upper bounded by a bias-variance
decomposition that also involves the similarity measure
$\simmeasgeneric$.

We begin by recalling the definition of the NW estimator, focusing
here on the version in which the underlying kernel is uniform over a
ball of a given bandwidth $\bandwidth_\numobs > 0$.  In particular,
define the set
\begin{equation*}
\goodset{\numobs} \defn \bigcup_{i=1}^\numobs
\ball{X_i}{\bandwidth_\numobs},
\end{equation*}
corresponding to the set of points in $\xspace$ within distance
$\bandwidth_\numobs$ of the observed covariates.  In terms of this
set, the \emph{Nadaraya-Watson estimator} $\fhatNW$ takes the form
\begin{equation}
\label{EqnNW}  
\fhat(x) \defn \begin{cases} \dfrac{\sum_{i=1}^\numobs Y_i \1\{X_i \in
    \ball{x}{\bandwidth_n}\}}{\sum_{i=1}^\numobs \1\{ X_i \in
    \ball{x}{\bandwidth_n}\}} & \mbox{for}~x \in \goodset{\numobs}
  \\ 0 & \mbox{otherwise.}
\end{cases}
\end{equation}
Our first main result provides an upper bound on the MSE of the NW
estimator under covariate shift; this bound exhibits the significance
of the similarity measure~\eqref{eqn:similarity-measure}.  It involves
the distribution $\midpointdist_\numobs \defn
\frac{\numsource}{\numobs} P + \frac{\numtarget}{\numobs} Q$, which is
a convex combination of the source and target distributions weighted
by their respective fractions of samples.
\btheo
\label{thm:upper-bound}
Suppose that Assumptions~\ref{asmp:holder}, \ref{asmp:covshift},
and~\ref{asmp:noise} hold.  For any $\bandwidth_\numobs > 0$, the
Nadaraya-Watson estimator $\fhatNW$ with bandwidth
$\bandwidth_\numobs$ has MSE bounded as
\begin{equation}
  \label{EqnNWBound}
  \E \big\|\fhatNW - \fstar\big\|_{L^2(\targetdist)}^2 \leq
  \upperconst \Big\{ L^2 \bandwidth_\numobs^{2\beta} +
  \frac{\|f^\star\|_\infty^2 + \noisebound^2}{\numobs}
  \similarity{\bandwidth_\numobs}{\midpointdist_n}{\targetdist}
  \Big\},
\end{equation}
where $\upperconst > 0$ is a numerical constant.
\etheo
\noindent See Section~\ref{sec:proof-upper} for a proof of this
result.

Note that the bound~\eqref{EqnNWBound} exhibits a type of
bias-variance trade-off, one that controls the optimal choice of
bandwidth $\bandwidth_\numobs$.  The quantity $\bandwidth_\numobs^{2
  \beta}$ in the first term is familiar from the classical analysis of
the NW estimator; it corresponds to the bias induced by smoothing over
balls of radius $\bandwidth_\numobs$, and hence is an increasing
function of bandwidth.  In the second term, the bandwidth appears in
the similarity measure
$\similarity{\bandwidth_\numobs}{\midpointdist_n}{\targetdist}$, which
is a non-increasing function of the bandwidth.  The optimal choice of
bandwidth arises from optimizing this tradeoff; note that it depends
on the pair $(P, Q)$, as well as the sample sizes $(\numobs_P,
\numobs_Q)$, via the similarity measure applied to the convex
combination $\mu_\numobs$ and $Q$.

\paragraph{No covariate shift:}  As a sanity check,
it is worth checking that the bound~\eqref{EqnNWBound} recovers known
results in the case of no covariate shift ($P = Q$ and hence
$\mu_\numobs = Q$). As a concrete example, if $Q$ is uniform on the
hypercube $[0,1]^k$, it can be verified that
$\similarity{\bandwidth}{\targetdist}{\targetdist} \asymp
\bandwidth^{-k}$ as $\bandwidth \rightarrow 0^+$.  (See
Example~\ref{ex:no-cov-shift} in the sequel for a more general
calculation that implies this fact.)  Thus, if we track only the
sample size, the optimal bandwidth is given by $\bandwidth^*_\numobs =
\numobs^{- \frac{1}{2 \beta + k}}$, and with this choice, the
bound~\eqref{EqnNWBound} implies that the NW estimator has MSE bounded
as $\numobs^{-\frac{2 \beta}{2 \beta + k}}$.  Thus, we recover the
classical and known results in this special case.  As we will see,
more interesting tradeoffs arise in the presence of covariate shift,
so that $\mu_\numobs \neq Q$.

\subsection{Consequences for $\alpha$-families of source-target pairs}
\label{SecAlphaFamily}

In order to better understand the bias-variance tradeoff in the
bound~\eqref{EqnNWBound} in the presence of covariate shift, it is
helpful to derive some explicit consequences of
Theorem~\ref{thm:upper-bound} for a particular function class $\cF$,
along with certain families of source-target pairs $(P, Q)$.  The
latter families are indexed by a parameter $\alpha > 0$ that controls
the amount of covariate shift; accordingly, we refer to them as
$\alpha$-families.

So as to simplify our presentation, we assume that $\xspace$ is the
unit interval $[0, 1]$.  For a given pair $\beta \in (0,1]$ and $L >
  0$, consider the class of regression functions
\begin{equation*}
\cF(\beta, L) = \Big\{\,f \colon [0, 1] \to \R \mid |f(x) - f(x')|
\leq L |x - x'|^\beta, \mbox{for all}~x, x' \in \xspace,~ f(0) = 0
\,\Big\}.
\end{equation*}
This is a special case of $\beta$-H\"{o}lder continuous functions when
the underlying metric space is the unit interval $[0,1]$ equipped with
the absolute value norm.  The additional constraint $f(0) = 0$ ensures
that this class has finite metric entropy. \\

\noindent Next we introduce some interesting families of source-target
pairs.

\paragraph{$\alpha$-families of $(P, Q)$ pairs:}
For a given parameter $\alpha \geq 1$ and radius $C \geq 1$, we define
the set of source-target pairs\footnote{Note that the restriction of
the supremum to $\bandwidth \in [0,1]$ is necessary, as
$\similarity{\bandwidth}{\sourcedist}{\targetdist} = 1$ for all
$\bandwidth \geq 1$. Note also that since
$\similarity{1}{\sourcedist}{\targetdist} = 1$, one necessarily has $C
\geq 1$.}
\begin{subequations}
  \begin{align}
\label{EqnAlphaBig}    
  \lowerboundclassbig{\alpha}{C} &\defn \Big\{\, (\sourcedist,
  \targetdist) \mid \sup_{0 < h \leq 1} h^{\alpha}
  \similarity{\bandwidth}{\sourcedist}{\targetdist} \leq C \,\Big\}.
\end{align}
In words, these are source target pairs for which the growth of the
similarity as $\bandwidth \rightarrow 0^+$ is at most
$\bandwidth^{-\alpha}$.  In the case $\alpha \in (0, 1]$, we define
  the related set
  \begin{align}
\label{EqnAlphaSmall}        
  \lowerboundclasssmall{\alpha}{C}&\defn \Big\{\, (\sourcedist,
  \targetdist) \mid \sup_{0 < h \leq \Delta} h^{\alpha}
  \similarity{\bandwidth}{\sourcedist}{\targetdist} \leq C, \sup_{0 <
    h \leq 1} \similarity{\bandwidth}{\targetdist}{\targetdist} \leq C
  \Big\},
\end{align}
\end{subequations}
where the additional condition is added to address the fact that even
without covariate shift, the rate $n^{-2\beta/(2\beta + 1)}$ is
unimprovable for some distributions~\cite{Sto82}.  Taking into account
the first part of the next corollary, it is necessary to impose some
condition on the target distribution in order to obtain significantly
faster rates such as $n^{-\tfrac{2\beta}{2\beta + \alpha}}$, when
$\alpha < 1$.
\bcor
\label{cor:consequences-upper}
Suppose that $\sigma \geq L$, and that Assumptions~\ref{asmp:covshift}
and~\ref{asmp:noise} hold.  Then there exists a constant $\upperconst'
> 0$, independent of $\numobs, \numsource, \numtarget, \sigma^2$, and
an integer $\numupper \defn \numupper(\sigma, \beta, L, \alpha, C)$
such that, provided that $\max\{\numsource, \numtarget\} \geq
\numupper$:
\begin{enumerate}[label={\upshape(\alph*)}]
\item For $\alpha \geq 1$ and $C \geq 1$, we have
\begin{subequations}
  \begin{equation}
\label{EqnUpperAlphaLarge}    
\inf_{\hat f} \sup_{\fstar \in \holderclass} \E \big\|\hat f -
\fstar\big\|_{L^2(\targetdist)}^2 \leq \upperconst' \Big \{
\big(\frac{\numsource }{\sigma^2} \big)^ \frac{2\beta + 1}{2\beta +
  \alpha} + \big(\frac{\numtarget} {\sigma^2} \big) \Big
\}^{\!-\frac{2\beta}{2\beta + 1}} ~~ \mbox{for any $(P, Q) \in
  \lowerboundclassbig{\alpha}{C}$.}
\end{equation}
\label{cor:alpha-big}
\item For $\alpha \in (0,1]$ and $C \geq 1$, we have
  \begin{equation}
\label{EqnUpperAlphaSmall}
\inf_{\hat f} \sup_{\fstar \in \holderclass} \E \big\|\hat f - \fstar
\big\|_{L^2(\targetdist)}^2 \leq \upperconst' \Big \{
\big(\frac{\numsource}{\sigma^2}\big)^{ \tfrac{2\beta}{2\beta +
    \alpha}} + \big(\frac{\numtarget}{\sigma^2}\big) \Big \}^{-1}
\quad \mbox{for any $(P, Q) \in \lowerboundclasssmall{\alpha}{C}$.}
\end{equation}
\end{subequations}
\label{cor:alpha-small}
\end{enumerate}
\ecor
\noindent See Section~\ref{sec:corollary-proof} for a proof of this
corollary. \\

Let us discuss the bound~\eqref{EqnUpperAlphaLarge} to gain some
intuition. The special case of no covariate shift can be captured by
setting $\numobs_P = 0$ and $\numobs_Q > 0$, and we recover the
familiar $\numobs^{- \frac{2 \beta}{2 \beta + k}}$ rate previously
discussed.  At the other extreme, suppose that $\numobs_Q = 0$ so that
all of our samples are from the shifted distribution (i.e., $\numobs =
\numobs_P$); in this case, the MSE is bounded as
$(\sigma^2/\numobs)^{-\frac{2 \beta}{2 \beta + \alpha}}$.  As $\alpha$
increases, our set-up allows for more severe form of covariate shift,
and its deleterious effect is witnessed by the exponent $\frac{2
  \beta}{2 \beta + \alpha}$ shrinking towards zero.  Thus, the NW
estimator---with an appropriate choice of bandwidth---remains
consistent but with an arbitrarily slow rate as $\alpha$ diverges to
$+\infty$.

\medskip

There are many papers in the literature (e.g.,~\cite{SugEtAl12,Kpo17}) that
discuss the covariate shift problem when the likelihood ratio is
bounded---that is, when $\targetdist$ is absolutely continuous with
respect to $\sourcedist$ and $\sup_{x \in \xspace} \tfrac{\ud
  \targetdist}{\ud \sourcedist}(x) \leq b$ for some $b \geq 1$.  We
say that the pair $(P, Q)$ are $b$-bounded in this case.

\bex[Bounded likelihood ratio]
\label{ex:bounded-LR}
Suppose that $\xspace = [0, 1]^\dimension$ with the Euclidean metric,
and consider a pair $(P, Q)$ with $b$-bounded likelihood ratio.  In
this special case, our general theory yields bounds in terms of the
$b$-weighted \emph{effective sample size}
\begin{equation}
  \numeff(b) \defn \frac{\numsource}{b} + \numobs_Q.
\end{equation}
In particular, it follows from the proof of
Corollary~\ref{cor:consequences-upper} that in the regime $\sigma^2
\geq L^2$, we have the upper bound
\begin{equation*}
\E \big\|\hat f - \fstar\big\|_{L^2(\targetdist)}^2 \leq \upperconst'
\Big(\frac{\sigma^2}{\numeff(b)}\Big)^{\frac{2\beta}{2\beta +
    \dimension}},
\end{equation*}
provided that $\numeff(b)$ is large enough.  Consequently, the effect
of covariate shift with $b$-bounded pairs is to reduce $\numobs_P$ to
$\numobs_P/b$.  Again, we recover the standard rate
$(\tfrac{\sigma^2}{\numobs})^{\tfrac{2\beta}{2\beta + \dimension}}$ in
the case of no covariate shift (or equivalently, when $b = 1$).  This
recovers a known result and is minimax optimal.  \eex


\subsection{Matching lower bounds for $\alpha$-families}
\label{SecLower}

Thus far, we have seen that the similarity measure $\simmeasgeneric$
plays a central role in determining the estimation error of the NW
estimator under covariate shift.  However, this is just one of many
possible estimators in nonparametric regression.  Does this similarity
measure play a more fundamental role?  In this section, we answer this
question in the affirmative by proving minimax lower bounds for
covariate shift problems parameterized in terms of bounds on
$\simmeasgeneric$.   In order to do so, we consider the metric
space $\xspace = [0,1]$ equipped with the absolute value as the
metric.

The main result of this section provides lower bounds on the
mean-squared error of any estimator, when measured uniformly over
functions in the H\"{o}lder class $\holderclass$, along with
target-source pairs $(P, Q)$ belonging to the class
$\lowerboundclassbig{\alpha}{C}$ when $\alpha \geq 1$ and the class
$\lowerboundclasssmall{\alpha}{C}$ when $\alpha < 1$.
\btheo
\label{thm:lower-bound}
Suppose that Assumptions \ref{asmp:covshift} and \ref{asmp:noise}
hold. Then there is a constant $\lowerconst > 0$, independent of
$\numobs, \numsource, \numtarget, \sigma^2$, and an integer $\numlower
\defn \numlower (\sigma, L, C, \alpha , \beta)$ such that for all
sample sizes $\max\{\numsource, \numtarget\} \geq \numlower$:
\begin{enumerate}[label={\upshape(\alph*)}]
\item
  \label{thm:alpha-big}
For $\alpha > 1$ and $C \geq 1$, there is a pair of
  distributions $(P, Q) \in \lowerboundclassbig{\alpha}{C}$ such that
  \begin{subequations}
\begin{equation}
  \label{ineq:lower-bound-ineq}
\inf_{\hat f} \sup_{\fstar \in \holderclass} \E \big\|\hat f -
\fstar\big\|_{L^2(\targetdist)}^2 \geq \lowerconst \Big \{
\big(\frac{\numsource }{\sigma^2} \big)^\frac{2\beta + 1}{2\beta +
  \alpha}+ \big( \frac{\numtarget} {\sigma^2} \big) \Big
\}^{-\frac{2\beta}{2\beta + 1}}.
\end{equation}
\item
  \label{thm:alpha-small}  
For $\alpha \leq 1$ and $C \geq 1$, there is a pair of
  distributions $(P, Q) \in \lowerboundclasssmall{\alpha}{C}$ such
  that
  \begin{equation}
  \label{ineq:lower-bound-ineq-two}    
\inf_{\hat f} \sup_{\fstar \in \holderclass} \E \big\|\hat f -
\fstar\big\|_{(\targetdist)}^2 \geq \lowerconst \Big \{ \big(\frac{
  \numsource}{\sigma^2} \big)^{\tfrac{2\beta}{2\beta + \alpha}} +
\big(\frac{\numtarget}{\sigma^2}\big) \Big \}^{-1}.
  \end{equation}
  \end{subequations}
\end{enumerate}
\etheo
\noindent See Sections~\ref{sec:proof-lower-big}
and~\ref{sec:proof-lower-small} for the proof of this result. \\

These lower bounds should be compared to
Corollary~\ref{cor:consequences-upper}. This comparison shows that the
MSE bounds achieved by the NW estimator are actually optimal in the
minimax sense over families defined by the similarity measure
$\simmeasgeneric$.


\section{Properties of the similarity measure}
\label{SecProperties}
In the previous sections, we have seen that the similarity measure
$\simmeasgeneric$ controls both the behavior of the NW estimator, as
well as fundamental (minimax) risks applicable to any estimator.
Thus, it is natural to explore the similarity measure in some more
detail, and in particular to draw some connections to existing notions
in the literature.

\subsection{Controlling $\simmeasgeneric$ via covering numbers}

We start with a general way of controlling the similarity measure
$\simmeasgeneric$, which is based on the covering number of the metric
space $(\xspace, d)$.  In particular, for any $h > 0$, the
\emph{covering number} $N(h)$ is defined to be the smallest number of
balls of radius $h$ needed to cover the space $\xspace$.  See Chapter
5 in the book~\cite{Wai19} for more background.

\bpr [Covering number bounds for the similarity measure]
\label{prop:covering-numbers}
Suppose that $P, Q$ are two probability measures on the same metric
space $(\xspace, d)$.  Suppose that for some $\bandwidth > 0$, there
is a $\lambda > 0$ such that
\begin{equation}
\label{EqnCoverBound}  
  P(\ball{x}{\bandwidth}) \geq \lambda~Q(\ball{x}{\bandwidth}) \quad
  \mbox{for all $x \in \xspace$.}
\end{equation}
Then the similarity at scale $\bandwidth$ is upper bounded as
$\similarity{\bandwidth}{P}{Q} \leq N(\tfrac{\bandwidth}{2})/\lambda$.
\epr

\noindent See Section~\ref{sec:proof-covering} for the proof of this claim. 
\medskip

It is worth emphasizing that---due to the order of quantifiers
above---the quantity $\lambda > 0$ is allowed to depend on $\bandwidth
> 0$.  We exploit this fact in subsequent uses of the
bound~\eqref{EqnCoverBound}. \\

One straightforward application of
Proposition~\ref{prop:covering-numbers} is in bounding the similarity
measure when there is no covariate shift, as we now discuss.

\bex [No covariate shift]
\label{ex:no-cov-shift}
Suppose that we compute the similarity measure in the case $P = Q$;
intuitively, this models a scenario where there is no covariate shift.
In this case, we clearly may apply
Proposition~\ref{prop:covering-numbers} with $\lambda = 1$, which
reveals that $\similarity{\bandwidth}{P}{P} \leq N(\bandwidth/2)$. To
give one concrete bound, suppose that $\xspace \subset \R^\dimension$
is a compact set, with diameter $D$. Then---owing to standard bounds
on covering number~\cite[chap.~5]{Wai19}---we obtain
$\similarity{\bandwidth}{P}{P} \leq (1 +
\tfrac{2D}{\bandwidth})^\dimension$. Note that this bound holds for
any metric, so long as the diameter $D$ is computed with the same
metric as the balls in the definition of the similarity measure.  \eex

\noindent We give another application of
Proposition~\ref{prop:covering-numbers} in the following subsection.

\subsection{Comparison to previous notions of distribution mismatch}
\label{sec:comparison}
Next, we show how the mapping $\bandwidth\mapsto
\similarity{\bandwidth}{P}{Q}$ can be bounded naturally using
previously proposed notions of distribution mismatch for covariate
shift. Again, Proposition~\ref{prop:covering-numbers} plays a central
role.

\bex [Bounded likelihood ratio]
Suppose that $P, Q$ are such that $Q \ll P$ and the likelihood ratio
$\tfrac{\ud Q}{\ud P}(x) \leq b$, for all $x \in \xspace$. Then note that
by a simple integration argument $P(\ball{x}{\bandwidth}) \geq 
\tfrac{1}{b} Q(\ball{x}{\bandwidth})$.
Therefore, we conclude $\similarity{\bandwidth}{P}{Q} \leq b N(\bandwidth/2)$. 
\eex

As noted previously, our work was inspired by the transfer exponent
introduced by Kpotufe and Martinet~\cite{KpoMar21} in the context of
covariate shift for nonparametric regression.  It is worth comparing
these notions so as to understand in what sense the similarity measure
$\simmeasgeneric$ is a refinement of the transfer exponent.  In order
to simplify this discussion, we focus here on the special case
$\xspace = [0,1]$.

\begin{figure}[h]
\centering 
\begin{overpic}[scale=0.6, unit=1mm]{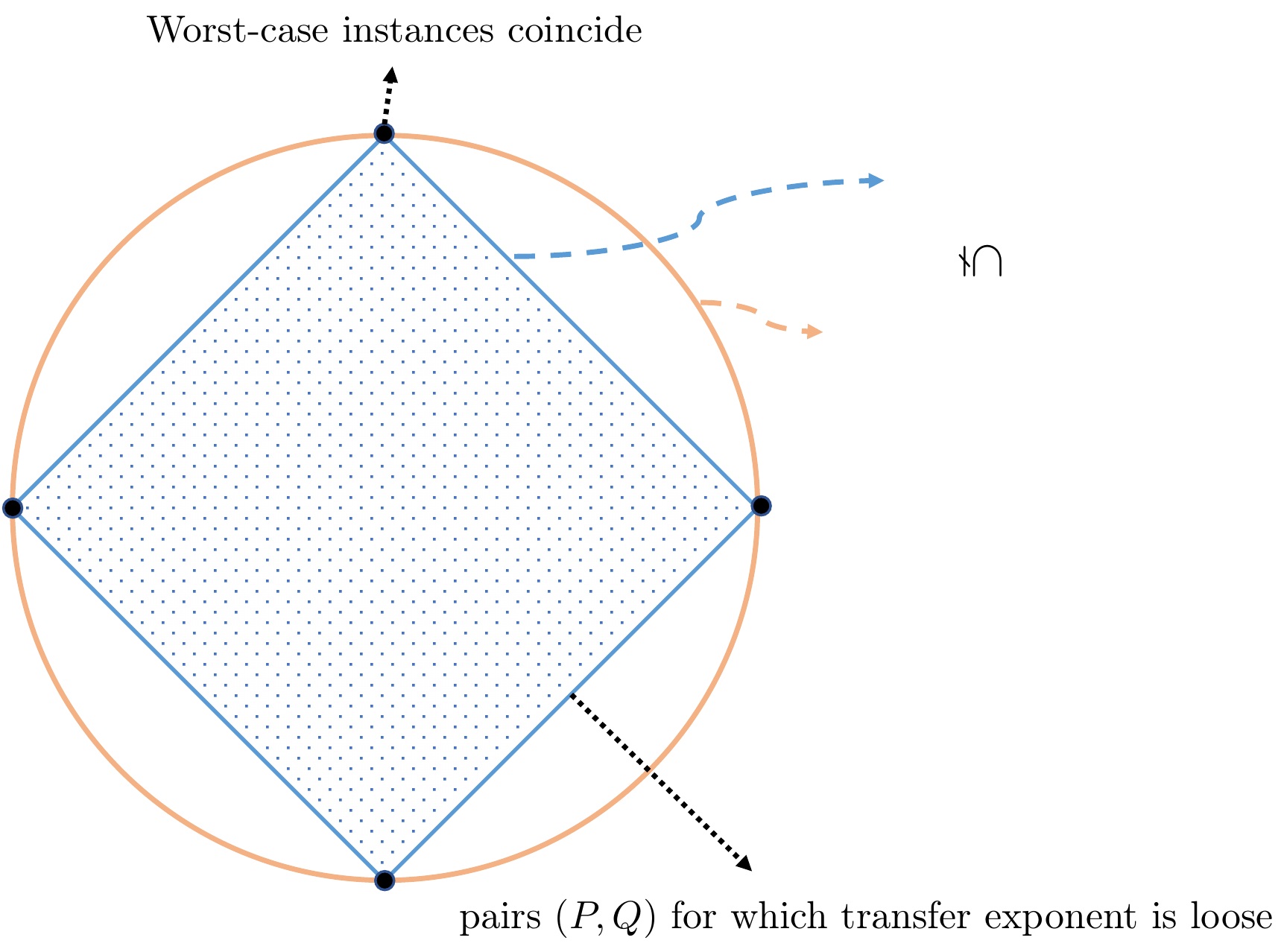}
\put(71, 59.3){$\cT(\gamma, \trad)$}
\put(67, 47.3){$\cD(\gamma + 1, \tfrac{2}{\trad})$}
\end{overpic}
  \caption{The yellow circle depicts the contour for the class
    $\mathcal{D}(\gamma + 1, \frac{2}{\trad})$, while the blue square plots
    the contour for the class $\mathcal{T}(\gamma, \trad)$. It can be seen
    from Lemma~\ref{lem:conversion-lemma} and Example~5 that
    $\mathcal{T}(\gamma, \trad)$ is strict subset of $\mathcal{D}(\gamma + 1,
    \frac{2}{\trad})$. In addition, our lower bound shows that under
    covariate shift, the worst-case instances for both classes
    coincide with each other. However, there exist instances $(P,Q)$
    where the characterization using transfer exponent is
    intrinsically loose.
  \label{fig:transfer-exponent}}
\end{figure}

\noindent We begin by providing the definition of transfer exponent:
\begin{definition}[Transfer exponent~\cite{KpoMar21}]
The distributions $(P,Q)$ have transfer exponent $\gamma \geq 0$ with
constant $\trad \in (0,1]$ if
\begin{equation*}
  P(\ball{x}{\bandwidth}) \geq \trad h^\gamma Q(\ball{x}{\bandwidth})
  \qquad \mbox{for all $x$ in the support of $Q$.}
\end{equation*}
\end{definition}
\noindent We denote by $\transClass(\gamma, \trad)$ the set of all
pairs $(P, Q)$ with this property.

It is natural to ask how the set $\transClass(\gamma, \trad)$ is
related to the $\alpha$-family previously defined in
equation~\eqref{EqnAlphaBig}.  The following result establishes an
inclusion:

\ble
\label{lem:conversion-lemma}
For $\xspace = [0, 1]$ and any $\gamma \geq 0$ and $\trad \in (0,1]$,
  we have the inclusion
  \begin{equation}
\label{EqnInclusion}    
\transClass(\gamma, \trad) \subset \lowerboundclassbig{\gamma +
  1}{\tfrac{2}{\trad}}.
\end{equation}
\ele
\noindent The proof of this inclusion is given in
Section~\ref{sec:proof-conversion}.  At a high level, it exploits
Proposition~\ref{prop:covering-numbers} to show that for any $(P,Q)
\in \transClass(\gamma, \trad)$, we have the bound
$\similarity{\bandwidth}{P}{Q} \leq \tfrac{1}{\trad \bandwidth^\gamma}
N(\bandwidth/2)$. \\

From the inclusion~\eqref{EqnInclusion}, it follows that any covariate
shift instance $(P,Q)$ with finite transfer exponent $\gamma \geq 0$
belongs to an $\alpha$-similarity family with $\alpha = \gamma + 1$.
In fact, following a proof similar to that of
Theorem~\ref{thm:lower-bound}, we can show that for $\gamma \geq 0$,
there is pair $(P, Q)$ in the class $\transClass(\gamma, \trad)$ such
that the minimax risk for $\beta$-H\"{o}lder-continous functions
scales as $\numobs_{P}^{-\frac{2 \beta}{ 2\beta + \gamma +1}}$.  Note
that this risk bound coincides with the minimax risk associated with
the class $\lowerboundclassbig{\gamma + 1}{\tfrac{2}{\trad}}$.  In
other words, from a \emph{worst case} point of view, the source-target
class $\transClass(\gamma, \trad)$ is equally as hard as the class
$\lowerboundclassbig{\gamma + 1}{\tfrac{2}{\trad}}$ for nonparametric
regression under covariate shift.  However, this worst case
equivalence does not capture the full picture: there are many
covariate shift families for which the transfer exponent provides an
overly conservative prediction, and so does not capture the
fundamental difficulty of the problem.  Let us consider a concrete
example to illustrate.

\bex [Separation between transfer exponent and $\simmeasgeneric$] Let
the target distribution $Q$ be a uniform distribution on the interval
$[0,1]$, and for some $\polyexp \geq 1$, suppose that the source
distribution $P$ has density $p(x) = (\polyexp +1) x^\polyexp$ for $x
\in [0,1]$. With these definitions, it can be verified that $(P,Q) \in
\transClass(\polyexp, \trad)$ for some constant $\trad \in (0,1]$, and
  moreover, that the quantity $\polyexp$ is the \emph{smallest
  possible} transfer exponent for this pair. In contrast, another
  direct computation shows that the pair $(P,Q)$ belongs to the class
  $\lowerboundclassbig{\polyexp }{C'}$ for some constant $C' > 0$.
  These two inclusions establish a separation between the rates
  predicted by the transfer exponent and the similarity
  $\simmeasgeneric$.  Indeed, as shown by our theory, the difficulty
  of estimation over $\lowerboundclassbig{\polyexp }{C'}$ is smaller
  than that prescribed by $\transClass(\polyexp, \trad)$. Indeed, if
  one observe $\numobs$ samples from the source distribution, the
  worst-case rate indicated by the computation from the transfer
  exponent is $\numobs^{-\frac{2\beta}{2\beta + \polyexp + 1}}$,
  whereas the rate guaranteed by the similarity measure
  $\simmeasgeneric$ is $\numobs^{-\frac{2\beta}{2\beta + \polyexp}}$.  As
  an explicit example, Lipschitz functions ($\beta = 1$) and $\polyexp
  = 1$, we obtain the slower rate $\numobs^{-1/2}$ versus the faster
  rate $\numobs^{-2/3}$, so that the ratio between the two rates
  diverges as $\numobs^{1/6}$ as the sample size grows.  \eex

\noindent See also Figure~\ref{fig:transfer-exponent} for an
illustration of the connections and differences between the similarity
measure and the transfer exponent.


\section{Proofs}
\label{SecProofs}

We now turn to the proofs of the results stated in the previous
section.

\subsection{Proof of Theorem~\ref{thm:upper-bound}}
\label{sec:proof-upper}

Recall that the estimate $\fhat$ depends on the observations $\{(X_i,
Y_i)\}_{i=1}^\numobs$, and so should be understood as a random
function.  The core of the proof involves proving that, for each $x
\in \xspace$, we have
\begin{equation}
\label{ineq:final-pointwise-upper}
\E\Big[\big(\fhatNW(x) - \fstar(x)\big)^2\Big] \leq L^2
\bandwidth_\numobs^{2 \beta} + 
\frac{4 \sigma^2 + \|\fstar\|_\infty^2}{\numobs}
\frac{1}{\midpointdist_\numobs(\ball{x}{\bandwidth_\numobs})},
\end{equation}
where the expectation is taking over the observations $\{(X_i,
Y_i)\}_{i=1}^\numobs$.  Given this inequality, the
claim~\eqref{EqnNWBound} of Theorem~\ref{thm:upper-bound} follows,
since by Fubini's theorem, we can write
\begin{align*}
\E \Big[\big\|\fhatNW - \fstar\big\|_{L^2(\targetdist)}^2\Big] =
\int_\xspace \E\Big[\big(\fhatNW(x) - \fstar(x)\big)^2\Big] \, \ud
Q(x).
\end{align*}
Applying inequality~\eqref{ineq:final-pointwise-upper} and recalling
the definition of the similarity measure yields the
claim~\eqref{EqnNWBound}.

We now focus on establishing the
bound~\eqref{ineq:final-pointwise-upper}.  Our proof makes use of the
conditional expectation of $\fhat$ given the covariates
\begin{equation*}
\fbar(x) \defn \E[\fhatNW(x) \mid X_1, \dots, X_\numobs],
\quad \mbox{for any}~x \in \xspace.
\end{equation*}
To be explicit, the expectation is taken over $Y_i \mid X_i$, $i = 1,
\ldots, \numobs$.  With this definition, our first result provides a
bound on the conditional bias and variance.
\ble
\label{lem:conditional-bounds}
For each $x \in \xspace$ almost surely, the Nadaraya-Watson estimator
$\fhatNW$ satisfies the bounds
\begin{subequations}
  \begin{align}
  \label{EqnCondBiasBd}    
    (\fbar(x) - \fstar(x) \big)^2 & \leq \|\fstar\|_\infty^2 \1\{x
  \not \in \goodset{\numobs}\} + L^2 \bandwidth_\numobs^{2\beta} \1\{x
  \in \goodset{\numobs}\} \qquad \mbox{and} \\
\label{EqnCondVarBd}
\E[(\fbar(x) - \fhatNW(x))^2 \mid X_1, \dots, X_\numobs] & \leq
\tfrac{\noisebound^2}{\sum_{i=1}^\numobs \1 \{X_i \in
  \ball{x}{\bandwidth_\numobs}\}} \1\{x \in \goodset{\numobs}\}.
  \end{align}
\end{subequations}
\ele
\noindent We prove this auxiliary claim at the end of this section. \\

Taking the results of Lemma~\ref{lem:conditional-bounds} as given, we
continue our proof of the bound~\eqref{ineq:final-pointwise-upper}.
For any fixed $x \in \xspace$, a conditioning argument yields
  \begin{equation*}  
  \E\big[(\fhatNW(x) - \fstar(x))^2\big] = \E\big[(\fbar(x) -
    \fstar(x) \big)^2\big] + \E\Big[ \E[(\fbar(x) - \fhatNW(x))^2 \mid
      X_1, \dots, X_\numobs] \Big].
  \end{equation*}
By applying the bounds~\eqref{EqnCondBiasBd} and~\eqref{EqnCondVarBd}
to the two terms above, respectively, we arrive at the upper bound
\mbox{$\E\big[(\fhatNW(x) - \fstar(x))^2\big] \leq \Term_1 +
  \Term_2$,} where
\begin{equation*}
\Term_1 \defn \|\fstar\|_\infty^2 \E[\1\{x \not \in
  \goodset{\numobs}\}] + L^2 \bandwidth_\numobs^{2\beta}, \quad
\mbox{and} \quad \Term_2 \defn \E \Big[\1\{x \in \goodset{\numobs}\}
  \tfrac{\noisebound^2}{\sum_{i=1}^\numobs \1\{X_i \in
    \ball{x}{\bandwidth_\numobs}\}}\Big].
\end{equation*}
We bound each of these terms in turn.

\paragraph{Bounding $\Term_1$:}
By definition, the set $\goodset{\numobs}$ involves $\numobs$
independent random variables, so that for any $x \in \xspace$, we have
\begin{align}
\label{ineq:prob-not-near-data}  
 \E\big[\1\{x \notin \goodset{\numobs}\}\big] = \Big(1 -
 \sourcedist\big(\ball{x}{\bandwidth_\numobs}\big)\Big)^\numsource
 \Big(1 -
 \targetdist\big(\ball{x}{\bandwidth_\numobs}\big)\Big)^\numtarget &
 \stackrel{{\rm (i)}}{\leq} \frac{1}{\numobs \,
   \midpointdist_\numobs(\ball{x}{\bandwidth_\numobs})},
\end{align}
where step (i) follows from the elementary inequality \mbox{$(1-p)^n
  (1-q)^m \leq \exp(-(np + mq)) \leq \tfrac{1}{np + mq}$,} valid for
$p, q \in (0, 1)$ and nonnegative integers $n, m$.  Consequently,
the first term is upper bounded as
\begin{subequations}
\begin{align}
\label{EqnTermOneUpper}
\Term_1 & \leq \|\fstar\|_\infty^2 \frac{1}{\numobs \,
  \midpointdist_\numobs(\ball{x}{\bandwidth_\numobs})} + L^2
\bandwidth_\numobs^{2\beta}.
\end{align}

\paragraph{Bounding $\Term_2$:}  For a fixed $x \in \xspace$, and for each $i = 1, \ldots, \numobs$, define the Bernoulli random variable
$Z_i = \1[X_i \in \ball{x}{\bandwidth_\numobs}] \in \{0,1 \}$, along
with the binomial random variables $U = \sum_{i=1}^{n_P} Z_i$ and $V =
\sum_{i = n_P + 1}^n Z_i$.  With these definitions, we can write
\begin{align*}
\sum_{i=1}^\numobs \1 \{ X_i \in \ball{x}{\bandwidth_\numobs} \} = U +
V, \quad \mbox{and} \quad \1 \{x \in \goodset{\numobs}\} = \1 \big\{ U
+ V > 0 \big \}.
\end{align*}
Consequently, by an elementary bound for binomial random variables
(see Lemma~\ref{lem:binomial}), it follows that
\begin{equation}
\label{EqnTermTwoUpper}
\Term_2 = \E \Big[\1 \{ U + V > 0 \} \frac{1}{U + V} \Big] \leq
\frac{4}{\numobs \,
  \midpointdist_\numobs(\ball{x}{\bandwidth_\numobs})}.
\end{equation}
\end{subequations}
Combining inequalities~\eqref{EqnTermOneUpper}
and~\eqref{EqnTermTwoUpper} yields the
claim~\eqref{ineq:final-pointwise-upper}. \\

\noindent The only remaining detail is to prove the auxiliary lemma
used in the proof.

\begin{proof}[Proof of Lemma~\ref{lem:conditional-bounds}]
Recall that by definition, we have
\begin{equation*}
\fbar(x) = \begin{cases} \dfrac{\sum_{i=1}^\numobs \fstar(X_i) \1\{X_i
    \in \ball{x}{\bandwidth_\numobs}\}}{\sum_{i=1}^\numobs \1\{ X_i
    \in \ball{x}{\bandwidth_\numobs}\}} & x \in \goodset{\numobs} \\ 0
  & x\notin \goodset{\numobs}
    \end{cases}
  \end{equation*}

\paragraph{Proof of the bound~\eqref{EqnCondBiasBd}:}
By a direct expansion, we have
  \begin{align*}
    \big(\fbar(x) - \fstar(x) \big)^2 \1\{x \in \goodset{\numobs}\} 
    &= \Big(\frac{\sum_{i=1}^\numobs (\fstar(x) - \fstar(X_i)) 
    \1\{X_i \in \ball{x}{\bandwidth_\numobs}\} }{\sum_{i=1}^\numobs 
    \1\{X_i \in \ball{x}{\bandwidth_\numobs}\}}\Big)^2
    \1\{x \in \goodset{\numobs}\} \\
    &\stackrel{{\rm (i)}}{\leq}
    \frac{\sum_{i=1}^\numobs (\fstar(x) - \fstar(X_i))^2 
    \1\{X_i \in \ball{x}{\bandwidth_\numobs}\}}
         {\sum_{i=1}^\numobs \1\{X_i \in \ball{x}{\bandwidth_\numobs}\}} 
         \1\{x \in \goodset{\numobs}\} \\
   &\stackrel{{\rm (ii)}}{\leq}
    L^2 \bandwidth_\numobs^{2\beta} 
    \1\{x \in \goodset{\numobs}\},
  \end{align*}
where step (i) follows from Jensen's inequality; and step (ii) makes
use of Assumption~\ref{asmp:holder}.  The bound~\eqref{EqnCondBiasBd}
is an immediate consequence.

\paragraph{Proof of the bound~\eqref{EqnCondVarBd}:}
In order to prove this claim, note that by independence among $\{(X_i,
\noise_i)\}_{i=1}^\numobs$,
  \begin{align*}
  \E[(\fbar(x) - \fhatNW(x))^2 \mid X_1, \dots, X_\numobs] &=
  \sum_{i=1}^\numobs
  \E[\noise_i^2 \mid X_i]
  \big(
  \tfrac{ \1\{X_i \in \ball{x}{\bandwidth_\numobs}\}}{ 
  \sum_{i=1}^\numobs \1\{ X_i \in \ball{x}{\bandwidth_\numobs} \}}
  \big)^2
  \1\{x \in \goodset{\numobs}\} \\
        &\stackrel{{\rm (iii)}}{\leq} \noisebound^2
        \sum_{i=1}^\numobs \big(\tfrac{\1\{X_i \in 
        \ball{x}{\bandwidth_\numobs}\}}
            {\sum_{i=1}^\numobs \1\{X_i \in 
            \ball{x}{\bandwidth_\numobs}\}}\big)^2  
            \1\{x \in \goodset{\numobs}\} \\
       &= \frac{\noisebound^2}{\sum_{i=1}^\numobs 
       \1\{X_i \in \ball{x}{\bandwidth_\numobs}\}} 
       \1\{x \in \goodset{\numobs}\},
  \end{align*}
  which proves the claim. Here step (iii) is a consequence of
  Assumption~\ref{asmp:noise}.
\end{proof}


\subsection{Proof of Corollary~\ref{cor:consequences-upper}}
\label{sec:corollary-proof}

Fix some $\bandwidth \in (0, 1]$, and introduce the indicator variable
  $\caseindicator = \1\{\alpha \geq 1\}$.  We then have
\begin{align*}
\int_{\xspace} \frac{1}{\numsource P(\ball{x}{h}) + \numtarget
  Q(\ball{x}{h})} \, \ud Q(x) & \leq \min\Big\{\frac{1}{\numsource}
\similarity{\bandwidth}{P}{Q} ,
\frac{1}{\numtarget}\similarity{\bandwidth}{Q}{Q}\Big\}\\ &\leq 3^\caseindicator
C \min\Big\{ \frac{1}{\numsource \bandwidth^\alpha},
\frac{1}{\numtarget \bandwidth^\caseindicator}\Big\} \\ &\leq 2 \cdot
3^\caseindicator C
\frac{1}{\numsource \bandwidth^\alpha + \numtarget \bandwidth^\caseindicator}.
\end{align*}
The last inequality follows from \eqref{prop:covering-numbers}
and standard covering number bounds (note $\bandwidth \leq 1$).
Thus the final performance bound is
\begin{equation*}
2 \cdot 3^\caseindicator C L^2 
\Big\{\bandwidth ^{2\beta}
+  \frac{L^2 + \sigma^2}{\numsource \bandwidth^\alpha + \numtarget
\bandwidth^\caseindicator}
\Big\}.
\end{equation*}
We choose the bandwidth $\bandwidth^\star$ so as to trade off between
two terms in this risk bound; more precisely, we set
\begin{equation*}
\bandwidth^\star = \Big(\big(\frac{\numtarget}{L^2 + \sigma^2}\big) +
\big(\frac{\numsource}{L^2 +
  \sigma^2}\big)^{\frac{2\beta+\caseindicator}{2\beta +
    \alpha}}\Big)^{-\tfrac{1}{2\beta + \caseindicator}}
\end{equation*}
This choice is valid, since $\sigma^2 \geq L^2$ and $\max\{\numsource,
\numtarget\} \geq 4\sigma^2$ by assumption.  Substituting this choice
of bandwidth into the risk bound~\eqref{EqnNWBound} yields the claim.

\subsection{Proof of Theorem~\ref{thm:lower-bound}\ref{thm:alpha-big}}
\label{sec:proof-lower-big}

Before giving the complete proof, we outline the main steps involved.
\begin{enumerate}
\item 
  We first construct a hard instance 
  $(P, Q) \in \lowerboundclassbig{\alpha}{C}$.
  This instance is designed such that the integral 
  quantity $\similarity{\bandwidth}{P}{Q}$ must scale as $C h^{-\alpha}$. 
\item Then we select a family of hard regression functions contained within 
$\cF(\beta, L)$ that guarantees the worst-case expected error for 
our pair of distributions, $(P, Q)$. 
\item Finally, we
  apply Fano's method over this set of regression functions to show that 
  the expected error must scale as the righthand side of 
  inequality~\eqref{ineq:lower-bound-ineq}. 
\end{enumerate}

It is worth commenting on our proof strategy in relation to past work.
On one hand, in the case $\alpha \geq 1$, our construction of the
distributions $(P, Q)$ is adapted from the lower bound argument
introduced by Kpotufe and Martinet~\cite{KpoMar21}.  The technical
work involves constructing pairs of densities of $P, Q$, and
establishing their membership in the class
$\lowerboundclassbig{\alpha}{C}$.  As for the case $\alpha \in (0,1)$,
as stated in Theorem~\ref{thm:lower-bound}\ref{thm:alpha-small}, we
use a different construction of the distribution pair $(P, Q)$, one
that is new (to the best of our knowledge).  We combine these
constructions of ``hard'' source-target pairs, in particular by
packing the interval $[0, 1]$ with a variable number of small
intervals (e.g.,~\cite{Yu93,Tsy09,Wai19}).  By adapting the number of
intervals (and constructing a packing set of the function class
$\holderclass$ appropriately over these intervals), one can adapt the
hardness of the lower bound instance to change with the number of
samples. In this case, we are able to do this such that the hardness
scales appropriately with the critical parameters that govern the
final minimax lower bound: $\numsource, \numtarget, \noisebound,
\alpha, \beta$.  With this high-level overview in place, we now
proceed to the technical content of the proof.

\paragraph{Constructing ``hard'' source-target pairs:}
\label{sec:hard-dist}
For scalars $\supportlength, r \in (0, 1]$, define $M =
  \tfrac{\supportlength}{6r}$ along with the intervals
\begin{equation*}
I_j \defn (z_j - 3r, z_j + 3r], \quad \mbox{where} \quad z_j \defn 6jr
  - 3r, \qquad j = 1, \ldots, M.
\end{equation*}
We specify $P$ and $Q$ on each interval $I_j$ as follows:
\begin{figure}[t]
\centering
  \includegraphics[scale=0.6]{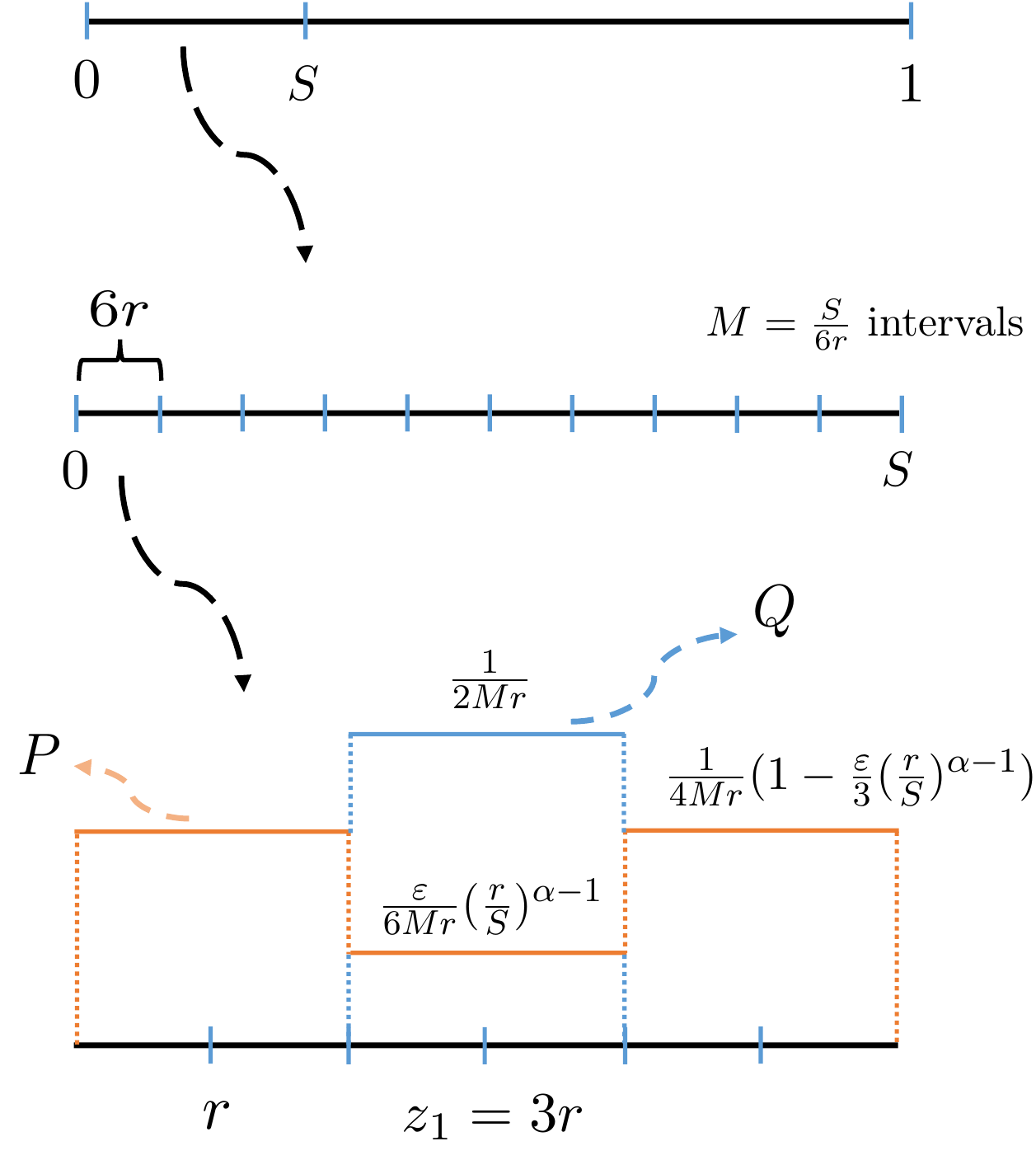}
  \caption{An illustration of the distributions $(P, Q)$ constructed as a hard pair in our lower bound.}
\end{figure}

\begin{table}[H]
  \centering
  \begin{tabular}{@{}lllll@{}}
    \toprule
    subinterval         & density of $P$                      & density of $Q$  &  &  \\ \midrule
    $(z_j-3r, z_j -r]$  & $\frac{1}{4Mr}(1 - \tfrac{\eps}{3}(\tfrac{r}{\supportlength})^{\alpha - 1})$ & 0               &  &  \\
      \addlinespace[0.5em]
      $(z_j-r, z_j +r]$   & $\frac{\eps}{6Mr} (\tfrac{r}{\supportlength})^{\alpha - 1}$      & $\frac{1}{2Mr}$ &  &  \\
        \addlinespace[0.5em]
        $(z_j +r, z_j +3r]$ & $\frac{1}{4Mr}(1 - \tfrac{\eps}{3}(\tfrac{r}{\supportlength})^{\alpha - 1})$ & 0               &  &  \\ \bottomrule
  \end{tabular}
  \caption{Specification of densities for lower bound pair of distributions $(P, Q)$ on the interval $I_j$.}
  \label{table:hard-distributions}
\end{table}
\noindent By construction, both $P$ and $Q$ assign probability 
$1/M$ to the entire interval $I_j$. The following proposition verifies that 
$(P,Q)$ lies in $\lowerboundclassbig{\alpha}{C}$ for proper choices 
of the $\eps$ and $S$.
\bpr
\label{prop:param-choice-lower}
Let $\alpha \geq 1$ and $C \geq 1$. Define $P$ and $Q$ as in 
Table~\ref{table:hard-distributions},
with the following choice of parameters $\eps, S$:
\begin{enumerate}[label={\upshape(\alph*)}]
\item if $C > 6$, set $\eps = 6/C$, and $S = 1/4$;
\item if $1 \leq C \leq 6$, set $\eps = 1$, and $S = \tfrac{1}{4}(C/6)^{1/\alpha}$.
\end{enumerate}
Then for any choice of $M, r > 0$ satisfying $S = 6Mr$, the pair
$(P,Q)$ lies in $\lowerboundclassbig{\alpha}{C}$.
\epr
\noindent See Section~\ref{sec:proof-param-choice} for the proof of
this claim.


\paragraph{Construction of ``hard'' regression functions.}
\label{sec:hard-fn}

Next we construct a packing of the function class of $\cF(\beta, L)$.
We do so by summing together scaled and shifted copies of base
function $\Psi \colon [-1, 1] \to \R$ that satisfies the boundary
conditions $\Psi(-1) = \Psi(1) = 0$, along with
\begin{subequations} 
\begin{align}
  \big|\Psi(x) - \Psi(y)\big| &\leq |x - y|^{\beta}, \quad \mbox{for
    all}~x, y \in [-1, 1],\quad \mbox{and},
  \label{ineq:psi-continuity}\\
  \int_{-1}^1 \Psi^2(x) \, \ud x &\eqcolon \psiintegral^2 > 0.
  \label{ineq:psi-integral}
\end{align}
\end{subequations}
There are many possible choices of $\Psi$; see Chapter 2 in the
book~\cite{Tsy09} for details.  For our proof, we also require the
bound $\psiintegral^2 \leq 1/6$, so that we make the explicit choice
\begin{equation*}
  \Psi(x) \defn e^{-1/(1-x^2)} \1\{|x| \leq 1\}.
\end{equation*}
We now form a class of functions using sums of the form
\begin{equation*}
f_\randbin(x) \defn \sum_{j=1}^M \randbin_j \phi_j(x), \quad
\mbox{where} \quad \phi_j(x) \defn L r^\beta \Psi\Big( \frac{x -
  z_j}{r}\Big),
\end{equation*}
and $\randbin = (\randbin_1, \dots, \randbin_M) \in \{0, 1\}^M$ is a
Boolean sequence.  Our construction makes use of the Gilbert-Varshamov
lemma (e.g. \cite[Lemma~2.9]{Tsy09}), which for $M \geq 8$, guarantees
the existence of a subset $\GVset \subset \{0, 1\}^M$ of cardinality
at least $2^{M/8}$ such that
\begin{equation}
\label{EqnGV}
\|\randbin - \randbin'\|_1 \geq M/8 \qquad \mbox{for all distinct $b,
  b' \in \GVset$.}
\end{equation}

\ble
\label{lem:function-class-properties}
The function class $\cH \defn \Big\{\, f_\randbin \mid \randbin \in
\GVset \, \Big\}$ has the following properties:
\begin{enumerate}[label={\upshape(\alph*)}]
\item It is contained within the H\"{o}lder class---$\cH \subset
  \cF(\beta, L)$.
  \label{item:prop-holder}
\item Pairs of functions are well-separated: for each distinct $f, g
  \in \cH$, we have
  \begin{align*}
    \|f - g\|^2_{L^2(\targetdist)} \geq \tfrac{\psiintegral^2}{16} L^2
    r^{2\beta}.
  \end{align*}
\label{item:prop-sep}      
\item Its elements satisfy the following $L^2(\sourcedist)$ and
  $L^2(\targetdist)$ bounds:
  \begin{equation*}
  \|f\|_{L^2(\targetdist)}^2 \leq \frac{\psiintegral^2
    M}{2\supportlength} L^2 r^{2\beta + 1} \quad \mbox{and} \quad
  \|f\|_{L^2(\sourcedist)}^2 \leq \frac{\eps \psiintegral^{2} M}{6
    S^{\alpha}} L^2 r^{2\beta + \alpha},
  \end{equation*}
  for all $f \in \cH$. \label{item:prop-bound}
\end{enumerate}
\ele

\paragraph{Applying Fano's method.}
\label{sec:lb-proof}

We now combine the preceding constructions with a Fano argument to
complete the proof of the lower bound.  For any function $f \in \cH$,
let $\jdist_f$ be the distribution $\{(X_i, Y_i)\}_{i=1}^\numobs$ where
$(X,Y)$ pairs are related by our nonparametric regression
model~\eqref{EqnObsModel} with $f = \fstar$.  For proving our lower
bound, it suffices to consider Gaussian noise: in particular,
$\noise_i \simiid \Normal{0}{\noisebound^2}$ for $i = 1, \ldots,
\numobs$.  These variables satisfy Assumption~\ref{asmp:noise}.

With these choices, Kullback-Leibler divergence between any given pair
$(\jdist_f, \jdist_g)$ can be bounded as
\begin{align*}
\kl{\jdist_f}{\jdist_g} & = \frac{1}{2 \noisebound^2} \Big(\numsource
\|f - g\|_{L^2(\sourcedist)}^2 + \numtarget \|f -
g\|_{L^2(\targetdist)}^2\Big) \leq \frac{2}{\noisebound^2}
\Big(\numsource \max_{f \in \cH} \|f\|_{L^2(\sourcedist)}^2 +
\numtarget \max_{f \in \cH}{\|f\|_{L^2(\targetdist)}^2} \Big).
\end{align*}
Now applying part~\ref{item:prop-bound} of
Lemma~\ref{lem:function-class-properties} yields
\begin{align*}
  \kl{\jdist_f}{\jdist_g} &\leq M\psiintegral^2 \Big\{
  \numsource \frac{L^2}{3\sigma^2} \frac{\eps}{S^\alpha}
  r^{2\beta + \alpha}
  + 
  \numtarget \frac{L^2}{\sigma^2} \frac{1}{S}
 r^{2\beta + 1}
  \Big\} \\
  &\leq M \Big\{ 
  \frac{4^\alpha}{C} \frac{L^2}{\sigma^2} \numsource 
  r^{2\beta + \alpha}
  + 
  \frac{4^\alpha}{C} \frac{L^2}{\sigma^2} \numtarget 
  r^{2\beta + 1}
  \Big\}  
\end{align*}
The final inequality arises by using $\psiintegral^2 \leq 1/6$. 
Suppose we take 
\begin{equation*}
r = \Big(\big(64 \frac{4^\alpha}{C} \frac{L^2 \numsource }{\sigma^2} )^
\frac{2\beta + 1}{2\beta + \alpha}+ \big(64  \frac{4^\alpha }{C} 
\frac{L^2 \numtarget}
{\sigma^2} \big) \Big)^{-\frac{1}{2\beta + 1}}
\end{equation*}
Then for any distinct $f, g \in \cH$, we obtain
\begin{equation}
  \label{ineq:KL-upper-bd} 
  \kl{\jdist_f}{\jdist_g} \leq M/32.  
\end{equation}
By a standard reduction to hypothesis testing~\cite[chap.~15]{Wai19} along 
with part~\ref{item:prop-holder}, 
\begin{equation*}
\inf_{\hat f}
\sup_{f^\star \cF(\beta, L)}
\E \Big[\|\hat f - \fstar\|_{L^2(\targetdist)}^2\Big]
\geq
\frac{\min_{(f, g) \in \binom{\cH}{2}} \|f - g\|_{L^2(\targetdist)}^2}{4}
  \Big\{1 - \frac{\log 2 + \max_{(f, g) \in 
  \binom{\cH}{2}} \kl{\jdist_f}{\jdist_g}}{\log |\cH|} \Big\}
  \end{equation*}
Thus, after applying part~\ref{item:prop-sep} of 
Lemma~\ref{lem:function-class-properties}, 
we obtain
\begin{equation*}
\inf_{\hat f}
\sup_{f^\star \cF(\beta, L)}
\E \Big[\|\hat f - \fstar\|_{L^2(\targetdist)}^2\Big]
\geq
\frac{\psiintegral^2}{64} L^2 r^{2\beta} \Big(1 - \frac{8}{M} - \frac{1}
{4}\Big) 
\geq 
\frac{\psiintegral^2 L^2}{256} 
\Big(\big(64 \frac{4^\alpha}{C} \frac{L^2 \numsource }{\sigma^2} )^
\frac{2\beta + 1}{2\beta + \alpha}+ \big(64  \frac{4^\alpha }{C} 
\frac{L^2 \numtarget}
{64\sigma^2} \big) \Big)^{-\frac{2\beta}{2\beta + 1}},
\end{equation*}
provided that $M \geq 32$. Equivalently, $r \leq \supportlength / 192$.
It suffices that $r \leq \tfrac{1}{4608}$, this is ensured by having 
\begin{equation*}
\max\{\numsource, \numtarget\}
\geq \Big(72 \frac{\sigma^2}{L^2} \frac{C}{4^\alpha}\Big)^{2\beta + \alpha}. 
\end{equation*}

\subsubsection{Proof of Proposition~\ref{prop:param-choice-lower}}
\label{sec:proof-param-choice}

We will show that for a general choice of 
$\eps, S \in (0, 1]$,
the following holds: 
\begin{equation}
\label{ineq:hard-dist-comparison-ineq}
P\big(\ball{x}{\bandwidth}\big) \geq \frac{\eps}{3} 
\big(\frac{\bandwidth}{4 \supportlength}\big)^{\alpha - 1} 
Q\big(\ball{x}{\bandwidth}\big),
\quad \mbox{for all}~x \in \mathrm{supp}(Q),~\mbox{and any}~h > 0. 
\end{equation}
For the moment let us take this bound as given.
By Lemma~\ref{lem:conversion-lemma}, note that 
bound~\eqref{ineq:hard-dist-comparison-ineq}
implies that $(P, Q) \in \lowerboundclassbig{\alpha}{\cC(\eps, S)}$,
with $\cC(\eps, S) = \tfrac{6}{\eps} (4\supportlength)^{\alpha -1}$,
for any $\eps, S \in (0, 1]$. 
Note that the parameter choices given in the statement of the 
result ensure that $\eps, S \in (0, 1]$.
  When $C \geq 6$, we have
  $\cC(\eps, S) = 6 (C/6)^{1 - 1/\alpha} = C (6/C)^{1/\alpha} \leq 6 \leq C$.
  Otherwise $C \leq 6$ and $\cC(\eps, S) = C$. 
  Therefore, checking the two cases $C > 6$ and $C \leq 6$ verifies $\cC(\eps, S)
= C$ in both regimes, which furnishes the claim.

We now turn to establish bound~\eqref{ineq:hard-dist-comparison-ineq}.
Let $h > 0$.  First observe that the support of $Q$ is the disjoint
union of intervals $\cup_{j = 1}^M (z_j -r, z_j+r]$. Thus, fix $x$ in
  the support of $Q$, and let $z_j$ denote the center of the interval
  to which $x$ belongs. Suppose that $\bandwidth \in [0, 4 r]$, in
  which case, we have the inclusion $\ball{x}{\bandwidth} \subset
  I_j$, whence the lower bound
\begin{align}
  P(\ball{x}{\bandwidth}) & \geq P\big(\ball{x}{\bandwidth} \cap
  \ball{z_j}{r}\big) \nonumber \\
  & \stackrel{{\rm(i)}}{=} \frac{\eps}{3}
  \Big(\frac{r}{\supportlength}\Big)^{\alpha - 1}
  Q\big(\ball{x}{\bandwidth} \cap
  \ball{z_j}{r}\big)\nonumber \\
  & \stackrel{{\rm(ii)}}{\geq} \frac{\eps}{3}
  \Big(\frac{\bandwidth}{4\supportlength}\Big)^{\alpha - 1}
  Q\big(\ball{x}{\bandwidth} \cap
  \ball{z_j}{r}\big)\nonumber\\
  & \stackrel{{\rm(iii)}}{=} \frac{\eps}{3}
  \Big(\frac{\bandwidth}{4\supportlength}\Big)^{\alpha - 1}
  Q\big(\ball{x}{\bandwidth}\big)
  \label{ineq:h-small}
\end{align}
Above, step (i) follows from the construction of $P, Q$; step (ii)
follows from $\bandwidth \leq 4r$, whereas step (iii) follows since
$\ball{x}{\bandwidth} \subset I_j$ and $Q$ assigns no mass to the set
$I_j \setminus \ball{z_j}{r}$.

Otherwise, we may assume that $\bandwidth \in [4r, \supportlength]$,
in which case we have the inclusion $\ball{x}{\bandwidth} \supset
I_j$. Denote by $N \geq 1$ the number of intervals of the form $I_j$
that are included within $\ball{x}{\bandwidth}$.  Note that since
$\ball{x}{\bandwidth}$ is connected, it is always contained in at most
$N + 2$ intervals (by considering partial intervals on the left and
right). Thus,
\begin{equation}\label{ineq:h-large}
\frac{P(\ball{x}{\bandwidth})}{Q(\ball{x}{\bandwidth})} \stackrel{{\rm(iii)}}
{\geq} \frac{N \cdot P(I_j)}{(N+2) \cdot Q(I_j)} \stackrel{{\rm(iv)}} \geq 
\frac{1}{3}. 
\end{equation}
Here step (iii) follows since $\ball{x}{\bandwidth}$ is contained in a
collection of at most $(N + 2)$ intervals and contains at least $N$
intervals, and the intervals are disjoint and have the same mass under
both $P$ and $Q$.  On the other hand, step (iv) uses the equivalence
$P(I_j) = Q(I_j)$, along with the fact that the function $x \mapsto
\tfrac{x}{x+2}$ is increasing on the set $\{x \geq 1\}$.

Therefore, combining inequalities~\eqref{ineq:h-small}
and~\eqref{ineq:h-large}, we conclude that
\begin{equation*}
P(\ball{x}{\bandwidth}) \geq \frac{1}{3} \Big[\eps
  \Big(\frac{\bandwidth} {4\supportlength}\Big)^{\alpha - 1} \wedge 1
  \Big] Q(\ball{x}{\bandwidth}) \geq \frac{\eps}{3}
\big(\frac{\bandwidth}{4 \supportlength}\big)^{\alpha - 1} Q(
\ball{x}{\bandwidth})
\end{equation*}
for every $x$ in the support of $Q$, the final inequality follows
since $\alpha \geq 1$. Since $h > 0$ was arbitrary, this establishes
bound~\eqref{ineq:hard-dist-comparison-ineq} and completes the proof.

\subsubsection{Proof of Lemma~\ref{lem:function-class-properties}}

We prove each of the three parts in turn.

\paragraph{Proof of part (a):} Fix a
Boolean vector $\randbin \in \{0, 1\}^M$.  Note that the function
$\phi_j$ is supported on the interval $I_j$, which is disjoint from
any other interval $I_k, k \neq j$.  Since $\Psi$ satisfies the
continuity condition~\eqref{ineq:psi-continuity}, it follows that
$\phi_j$ is $(\beta, L)$-Hölder.  Finally, we have $f_\eps(0) = 0$ by
definition.  Taking these properties together, we have shown that
$f_\eps \in \cF(\beta, L)$, as required.

\paragraph{Proof of part~\ref{item:prop-sep}:}
For any distinct pair $\randbin, \randbin' \in \GVset$,
we have
\begin{align*}
  \int_0^1 (f_{\randbin}(x) - f_{\randbin'}(x))^2 \, \ud Q(x) &=
  \int_0^1 \Big(\sum_{j=1}^M (\randbin_j - \randbin'_j) \phi_j(x)
  \Big)^2 \, \ud Q(x) \\
  & \stackrel{{\rm (i)}}{=} \frac{1}{2 M r} \sum_{j=1}^M (\randbin_j -
  \randbin'_j)^2 \int_{z_j - 3r}^{z_j + 3r} \phi_j^2(x) \, \ud x \\
    & \stackrel{{\rm (ii)}}{=} \frac{\psiintegral^2}{2M} L^2
  r^{2\beta} \|\randbin - \randbin'\|_1 \\
    & \stackrel{{\rm (iii)}}{\geq} \frac{\psiintegral^2}{16} L^2
  r^{2\beta}.
  \end{align*}
Here step (i) follows from the definition of $Q$ along with the
disjointedness of the supports of $\phi_j$. Step (ii) follows from
equation~\eqref{ineq:psi-integral} and the fact that $\randbin,
\randbin' \in \GVset \subset \{0, 1\}^M$.  Finally, step (iii) follows
from the Gilbert-Varshamov separation~\eqref{EqnGV}.

\paragraph{Proof of part~\ref{item:prop-bound}:}
For any $b \in \GVset$, by following the calculations above, for $\mu
\in \{P, Q\}$, we have by symmetry
  \begin{equation*}
  \int_{0}^1 f_b^2(x) \, \ud \mu(x) = \sum_{j=1}^M \randbin_j^2
  \int_{I_j} \phi_j^2(x) \, \ud \mu(x) \leq M \int_{I_1} \phi_1^2(x)
  \, \ud \mu(x).
  \end{equation*}
Now observe that $\int_{0}^{6r} \phi_1^2(x) \, \ud Q(x) =
\tfrac{\psiintegral^2}{2M} L^2 r^{2\beta}$, and consequently,
$\|f_b\|_{L^2(\targetdist)}^2 \leq L^2 r^{2\beta} \psiintegral^2/2$.
Additionally, we can compute
\begin{equation*}
\int_0^{6r} \phi_1^2(x) \, \ud P(x) = \frac{\eps}{6 r M^\alpha}
\int_{2r}^{4r} \phi_1^2(x) \, \ud x =
\frac{\eps}{6\supportlength^\alpha} L^2 r^{2\beta + \alpha}
\psiintegral^2.
\end{equation*} 
Thus, we have established the upper bound
$\|f_b\|_{L^2(\sourcedist)}^2 \leq \eps L^2 r^{2\beta + \alpha - 1}/(6
S^{\alpha - 1})$.


\subsection{Proof of Theorem~\ref{thm:lower-bound}\ref{thm:alpha-small}}
\label{sec:proof-lower-small}

Given the inclusion $\lowerboundclasssmall{\alpha}{1} \subset
\lowerboundclasssmall{\alpha}{C}$, it suffices to prove a lower bound
for $C = 1$.

\paragraph{Construction of ``hard'' distributions.}
Let $Q = \delta_1$, and let $P_\alpha$ be the distribution supported
on $[0,1]$ with density $p_\alpha(x) \defn \alpha (1-x)^{\alpha - 1}
\1\{x \in [0, 1]\}$.  By construction, we then have
\begin{equation*}
  \similarity{h}{P_\alpha}{Q} = \frac{1}{P_\alpha(B(1, h))} =
  h^{-\alpha} \quad \mbox{for all $h \in (0,1]$,}
\end{equation*}
which implies that $(P_\alpha, Q) \in
\lowerboundclasssmall{\alpha}{1}$.  From herein, we adopt the
shorthand $P \defn P_\alpha$ so as to lighten notation.

\paragraph{Construction of two point alternative.}
If the regression function is $f$, we denote the resulting joint
distribution of $\{(X_i, Y_i)\}_{i=1}^\numobs$ by $\jdist_f$.  We
consider the two point alternatives $\{f_t, g\}$ with $g \equiv 0$ and
$f_t(x) \defn L(x-t)_+^\beta$.  The next result demonstrates the
validity of this choice:
\ble
\label{lem:holder-part-small}
For any $t \in [0, 1]$, the function $f_t$ belongs to $\cF(\beta, L)$.
\ele
\noindent See Section~\ref{sec:proof-holder-part-small} for the
proof. \\

Moreover, by straightforward calculations, we find that
$\|f_t\|_{L^2(Q)}^2 = L^2 (1-t)^{2 \beta}$, and
\begin{align*}
\|f_t\|_{L^2(P)}^2 &= L^2 \int_t^1 \alpha (1-x)^{\alpha - 1}
(x-t)^{2\beta}\, \ud x  \\ &\leq L^2 (1-t)^{2\beta} \int_{0}^{1-t} \alpha
s^{\alpha -1}\, \ud s = L^2 (1 - t)^{2\beta + \alpha}.
\end{align*}

\paragraph{Applying Le Cam's method.} We are now equipped to apply Le Cam's two point bound.  In particular, we have
\begin{equation*}
\inf_{\hat f} \sup_{f^\star \in \cF(\beta, L)} \E\Big[ \|\hat f -
  f^\star\|_{L^2(Q)}^2\Big] \geq \frac{L^2(1-t)^{2\beta}}{16}
\exp\big(-\kl{\jdist_{f_t}}{\jdist_g}\big)
\end{equation*}
By standard KL calculations (using $\Normal{0}{\sigma^2}$ noises)
\begin{equation*}
\kl{\jdist_{f_t}}{\jdist_g} = \frac{L^2}{2\sigma^2}\Big\{ \numsource
(1-t)^{2\beta + \alpha} + \numtarget (1-t)^{2\beta} \Big\}
\end{equation*}
Finally, we make the
\begin{equation*}
1-t = \bigg( \Big(\frac{L^2
  \numsource}{2\sigma^2}\Big)^{\tfrac{1}{2\beta + \alpha}} +
\Big(\frac{L^2
  \numtarget}{2\sigma^2}\Big)^{\frac{1}{2\beta}}\bigg)^{-1}
\end{equation*}
A little bit of algebra shows that this choice guarantees that
$\kl{\jdist_{f_t}}{\jdist_g} \leq 2$, which completes the proof.


\subsubsection{Proof of Lemma~\ref{lem:holder-part-small}}
\label{sec:proof-holder-part-small}

We begin by observing that $f_t(0) = 0$. Thus, in order to prove the
claim, it suffices to show that
\begin{equation*}
  f_t(y) - f_t(x) \leq L (y - x)^\beta \quad \mbox{for any pair $x, y$
    such that $0 \leq t < x < y \leq 1$.}
\end{equation*}
In order to prove this bound, consider an arbitrary point $x \in (t,
1)$, and define the function
\begin{equation*}
\phi_x(y) \defn L(y^\beta - x^\beta) - L(y-x)^\beta \qquad \mbox{for
  $y \in [x, 1]$.}
\end{equation*}
We can compute the derivative $\phi_x'(y) = L\beta( y^{\beta - 1} -
(y-x)^{\beta - 1})$.  Since $y \geq y - x > 0$ and $\beta \leq 1$, we
have $y^{\beta-1} \leq (y-x)^{\beta-1}$, and hence $\phi_x'(y) \leq
0$. Consequently, the function $\phi_x$ is non-increasing, and since
$y > x$, it follows that $\phi_x(y) \leq \phi_x(x) = 0$. Putting
together the pieces completes the proof.

\subsection{Proof of Proposition~\ref{prop:covering-numbers}}
\label{sec:proof-covering}

Starting with the assumed bound~\eqref{EqnCoverBound}, we have
\begin{align}
\label{EqnCivetCoffee}  
\int_\xspace \frac{1}{P(\ball{x}{h})} \, \ud Q(x) &\leq
\frac{1}{\lambda} \int_\xspace \frac{1}{Q(\ball{x}{\bandwidth})} \,
\ud Q(x).
\end{align}
By definition of the covering number $N \defn N(\bandwidth/2)$, there
is a collection $\{z^j\}_{j=1}^N$ such that the set $\xspace$ is
contained within the union $\bigcup_{j=1}^N
\ball{z^j}{\frac{\bandwidth}{2}}$.  This fact, combined with our
previous bound~\eqref{EqnCivetCoffee}, implies that
\begin{align}
  \label{EqnYeastCoffee}
\int_\xspace \frac{1}{P(\ball{x}{h})} \, \ud Q(x) &\leq
\frac{1}{\lambda} \sum_{j=1}^N \int_{\ball{z_j}{\bandwidth/2}}
\frac{1}{Q(\ball{x}{\bandwidth})}\, \ud Q(x).
\end{align}
Note by the triangle inequality, for each $j \in [N]$ and $x \in
\ball{z_j}{\bandwidth/2}$, we have $\ball{z_j}{\bandwidth/2} \subset
\ball{x}{\bandwidth}$.  This inclusion implies that
\begin{equation*}
  \int_{\ball{z_j}{\bandwidth/2}} \frac{1}{Q(\ball{x}{\bandwidth})}\,
  \ud Q(x) \leq \int_{\ball{z_j}{\bandwidth/2}}
  \frac{1}{Q(\ball{z_j}{\bandwidth/2})}\, \ud Q (x) =1,
\end{equation*}
for each $j \in [N]$. Combining this inequality with the
bound~\eqref{EqnYeastCoffee} yields the claim.

\subsection{Proof of Lemma~\ref{lem:conversion-lemma}}
\label{sec:proof-conversion}

By assumption, we have the upper bound
\begin{equation*}
\int_0^1 \frac{1}{P(\ball{x}{\bandwidth})} \, \ud Q(x) \leq
\frac{1}{\trad h^\gamma} \int_0^1 \frac{1}{Q(\ball{x}{\bandwidth})} \,
\ud Q(x)
\end{equation*}
Moreover, we can find a collection of $N \defn \ceil{1/\bandwidth}$
balls with centers $\{z_j\}_{j=1}^N$ of radius $\bandwidth/2$ that
cover the interval $[0, 1]$, whence
\begin{equation*}
\int_0^1 \frac{1}{Q(\ball{x}{\bandwidth})} \, \ud Q(x) \leq
\sum_{j=1}^N \int_{x \in \ball{z_j}{\bandwidth/2}}
\frac{1}{Q(\ball{x}{\bandwidth})} \, \ud Q(x) \leq N.
\end{equation*}
The final inequality follows from the inclusion $\ball{x}{\bandwidth}
\supset \ball{z_j}{\bandwidth/2}$.

Now define the function $g(t) \defn \ceil{t}/t$, and observe that
$g(t) \leq 2$ whenever $t \geq 1$. Consequently, we can write
\begin{equation*}
\bandwidth^{\gamma + 1} \similarity{\bandwidth}{P}{Q} \leq
\frac{1}{\trad} g(1/\bandwidth) \leq \frac{2}{\trad}, \quad \mbox{for
  any}~h \leq 1.
\end{equation*}
Passing to the supremum over $\bandwidth \in (0, 1]$ yields the claim.

  
\section{Discussion}
\label{SecDiscussion}

In this paper, we have studied the problem of covariate shift in the
context of nonparametric regression.  We have shown that a measure of
(dis)-similarity $\simmeasgeneric$ between the source and target
distributions, as defined in equation~\eqref{eqn:similarity-measure},
can be used to characterize how minimax risks change as the
source-target pair are varied.  In particular, we proved upper bounds
on the Nadaraya-Watson estimator over H\"{o}lder classes that are an
explicit function of the similarity $\simmeasgeneric$, and also
established matching lower bounds over classes constrained in terms of
the similarity.  We also discussed how the measure $\simmeasgeneric$
is related to other characterizations of covariate shift from past
work, including likelihood ratio bounds and transfer exponents.  Our
work shows that similarity measure $\simmeasgeneric$ provides a more
fine-grained characterization of how covariate shift changes the
difficulty of non-parametric regression.

Our work leaves open a number of open questions.  First, our lower
bounds for covariate shift (cf. Theorem~\ref{thm:lower-bound}) are
obtained within a global minimax framework, which involves worst-case
assessments over a certain function class.  These lower bounds match
our upper bound on the NW estimator
(cf. Theorem~\ref{thm:upper-bound}) for certain source-target pairs
$(P, Q)$.  But the upper bound actually depends explicitly on the
source-target pair.  Is this upper bound always optimal? Or are there
instances of covariate shift for which Nadaraya-Watson is suboptimal
for some H\"{o}lder continuous function?  In general, this question
appears non-trivial: even without the (interesting) complication of
covariate shift, there are few results that give
distribution-dependent results for nonparametric regression outside of
the uniform distribution and fixed-design problems.


\subsection*{Acknowledgements}

We would like to thank Samory Kpotufe for helpful email exchanges.
This work was partially supported by NSF-DMS grant 2015454, NSF-IIS
grant 1909365, NSF-FODSI grant 202350, and DOD-ONR Office of Naval
Research N00014-21-1-2842 to MJW.

\appendix

\section{Elementary bound for binomial variables}

In this section, we state and prove an elementary bound for binomial
random variables, used in the proof of Theorem~\ref{thm:upper-bound}.
\ble
\label{lem:binomial}
Let $n, m$ be positive integers and $p, q \in (0, 1)$.
Suppose that
$U \sim \Bin{n}{p}$ and $V \sim \Bin{m}{q}$.  Then
\begin{equation*}
\E\Big[ \frac{1}{U+V} \1\{U + V > 0\}\Big]
\leq \frac{4}{np + mq}.
\end{equation*}
\ele
\begin{proof}
We begin by observing that conditionally on the event $\{U + V > 0\}$,
we have the lower bound
\begin{equation*}
  U + V \geq \frac{U + V + 1}{2} \geq \frac{U+1}{2} \vee
  \frac{V+1}{2}.
\end{equation*}
These lower bounds allow us to write
\begin{equation*}
\E \Big[ \frac{1}{U+V} \1\{U + V > 0\}\Big] \leq \E \frac{2}{U+1}
\wedge \E \frac{2}{ V + 1} \leq \frac{2}{(np \vee mq)} \leq
\frac{4}{np + mq}.
\end{equation*}
Here the penultimate inequality is a consequence of known results for
binomial random variables~\cite[equation (3.4)]{ChaStr72}.
\end{proof}

\bibliographystyle{plain}
\bibliography{references}

\begin{thebibliography}{10}

\bibitem{BenBliEtAl10}
S.~Ben-David, J.~Blitzer, K.~Crammer, A.~Kulesza, F.~Pereira, and J.~W.
  Vaughan.
\newblock A theory of learning from different domains.
\newblock {\em Machine learning}, 79(1-2):151--175, 2010.

\bibitem{BenLuEtAl10}
S.~Ben-David, T.~Lu, T.~Luu, and D.~P{\'a}l.
\newblock Impossibility theorems for domain adaptation.
\newblock In {\em Proceedings of the Thirteenth International Conference on
  Artificial Intelligence and Statistics}, pages 129--136, 2010.

\bibitem{BliEtAl07}
J.~Blitzer, M.~Dredze, and F.~Pereira.
\newblock Biographies, {B}ollywood, boom-boxes and blenders: Domain adaptation
  for sentiment classification.
\newblock In {\em Proceedings of the 45th Annual Meeting of the Association of
  Computational Linguistics}, pages 440--447, Prague, Czech Republic, June
  2007. Association for Computational Linguistics.

\bibitem{ChaStr72}
M.~T. Chao and W.~E. Strawderman.
\newblock Negative moments of positive random variables.
\newblock {\em J. Amer. Statist. Assoc.}, 67(338):429--431, 1972.

\bibitem{CorEtAl19}
C.~Cortes, M.~Mohri, and A.~M. Medina.
\newblock Adaptation based on generalized discrepancy.
\newblock {\em J. Mach. Learn. Res.}, 20(1):1--30, 2019.

\bibitem{GerEtAl13}
P.~Germain, A.~Habrard, F.~Laviolette, and E.~Morvant.
\newblock A pac-bayesian approach for domain adaptation with specialization to
  linear classifiers.
\newblock In {\em International Conference on Machine Learning}, pages
  738--746, 2013.

\bibitem{Hal72}
G.~Hal\'{a}sz.
\newblock Statistical interpolation.
\newblock {\em Mat. Lapok}, 23:71--87 (1973), 1972.

\bibitem{HasEtAl13}
A.~Hassan, R.~Damper, and M.~Niranjan.
\newblock On acoustic emotion recognition: Compensating for covariate shift.
\newblock {\em IEEE Trans. Audio Speech Lang. Process.}, 21(7):1458--1468,
  2013.

\bibitem{IbrKha80}
I.~A. Ibragimov and R.~Z. Khas\cprime~minski\u{\i}.
\newblock Nonparametric regression estimation.
\newblock {\em Dokl. Akad. Nauk SSSR}, 252(4):780--784, 1980.

\bibitem{Kpo17}
S.~Kpotufe.
\newblock Lipschitz density-ratios, structured data, and data-driven tuning.
\newblock In {\em Artificial Intelligence and Statistics}, pages 1320--1328,
  2017.

\bibitem{KpoMar21}
S.~Kpotufe and G.~Martinet.
\newblock Marginal singularity and the benefits of labels in covariate-shift.
\newblock {\em Ann. Statist.}, 49(6):3299--3323, 2021.

\bibitem{LeiHuLee21}
Q.~Lei, W.~Hu, and J.~Lee.
\newblock Near-optimal linear regression under distribution shift.
\newblock In Marina Meila and Tong Zhang, editors, {\em Proceedings of the 38th
  International Conference on Machine Learning}, volume 139 of {\em Proc. Mach.
  Learn. Res.}, pages 6164--6174. PMLR, 18--24 Jul 2021.

\bibitem{LiEtAl2010}
Y.~Li, H.~Kambara, Y.~Koike, and M.~Sugiyama.
\newblock Application of covariate shift adaptation techniques in
  brain–computer interfaces.
\newblock {\em IEEE Trans. Biomed. Eng.}, 57(6):1318--1324, 2010.

\bibitem{ManEtAl09}
Y.~Mansour, M.~Mohri, and A.~Rostamizadeh.
\newblock Domain adaptation: Learning bounds and algorithms.
\newblock {\em arXiv preprint arXiv:0902.3430}, 2009.

\bibitem{ManMohEtAl09}
Y.~Mansour, M.~Mohri, and A.~Rostamizadeh.
\newblock Multiple source adaptation and the r{\'e}nyi divergence.
\newblock In {\em Proceedings of the Twenty-Fifth Conference on Uncertainty in
  Artificial Intelligence}, pages 367--374. AUAI Press, 2009.

\bibitem{MohEtAl12}
M.~Mohri and A.~M. Medina.
\newblock New analysis and algorithm for learning with drifting distributions.
\newblock In {\em International Conference on Algorithmic Learning Theory},
  pages 124--138. Springer, 2012.

\bibitem{MouMohEtAl20}
M.~Mousavi~Kalan, Z.~Fabian, S.~Avestimehr, and M.~Soltanolkotabi.
\newblock Minimax lower bounds for transfer learning with linear and one-hidden
  layer neural networks.
\newblock In H.~Larochelle, M.~Ranzato, R.~Hadsell, M.~F. Balcan, and H.~Lin,
  editors, {\em Advances in {N}eural {I}nformation {P}rocessing {S}ystems},
  volume~33, pages 1959--1969. Curran Associates, Inc., 2020.

\bibitem{Nad64}
E.~A. Nadaraya.
\newblock On estimating regression.
\newblock {\em Theory Probab. Appl.}, 9(1):141--142, 1964.

\bibitem{QuiEtAl09}
J.~Quionero-Candela, M.~Sugiyama, A.~Schwaighofer, and N.~D. Lawrence.
\newblock {\em Dataset Shift in Machine Learning}.
\newblock The MIT Press, 2009.

\bibitem{SaeEtAl2010}
K.~Saenko, B.~Kulis, M.~Fritz, and T.~Darrell.
\newblock Adapting visual category models to new domains.
\newblock In K.~Daniilidis, P.~Maragos, and N.~Paragios, editors, {\em Computer
  Vision -- ECCV 2010}, pages 213--226, Berlin, Heidelberg, 2010. Springer
  Berlin Heidelberg.

\bibitem{Shi00}
H.~Shimodaira.
\newblock Improving predictive inference under covariate shift by weighting the
  log-likelihood function.
\newblock {\em J. Statist. Plann. Inference}, 90(2):227--244, 2000.

\bibitem{Sto82}
C.~J. Stone.
\newblock Optimal global rates of convergence for nonparametric regression.
\newblock {\em Ann. Statist.}, 10(4):1040--1053, 1982.

\bibitem{SugNakEtAl08}
M.~Sugiyama, S.~Nakajima, H.~Kashima, P.~V. Buenau, and M.~Kawanabe.
\newblock Direct importance estimation with model selection and its application
  to covariate shift adaptation.
\newblock In {\em Advances in {N}eural {I}nformation {P}rocessing {S}ystems},
  pages 1433--1440, 2008.

\bibitem{SugEtAl12}
M.~Sugiyama, T.~Suzuki, and T.~Kanamori.
\newblock {\em Density ratio estimation in machine learning}.
\newblock Cambridge University Press, 2012.

\bibitem{Tsy09}
A.~B. Tsybakov.
\newblock {\em Introduction to nonparametric estimation}.
\newblock Springer Series in Statistics. Springer, New York, 2009.
\newblock Revised and extended from the 2004 French original, Translated by
  Vladimir Zaiats.

\bibitem{Wai19}
M.~J. Wainwright.
\newblock {\em High-dimensional statistics: A non-asymptotic viewpoint}.
\newblock Cambridge University Press, Cambridge, {U}{K}, 2019.

\bibitem{Wat64}
G.~S. Watson.
\newblock Smooth regression analysis.
\newblock {\em Sankhy{\=a}: The Indian Journal of Statistics, Series A}, pages
  359--372, 1964.

\bibitem{Yu93}
B.~Yu.
\newblock Assouad, {F}ano and {L}e {C}am.
\newblock In {\em Festschrift in Honor of {L}. {L}e {C}am on his 70th
  Birthday}. 1993.

\end{thebibliography}
\end{document}